\documentclass[11pt]{article}
\usepackage{amsmath}
\usepackage{amsfonts}
\usepackage{amssymb}

\newtheorem{theorem}{Theorem}[section]
\newtheorem{lemma}[theorem]{Lemma}
\newtheorem{proposition}[theorem]{Proposition}
\newtheorem{definition}{Definition}[section]

\newtheorem{hypothesis}[theorem]{Hypothesis}

\newtheorem{remark}[theorem]{Remark}
\newtheorem{corollary}[theorem]{Corollary}

 \makeatletter
  \@addtoreset{equation}{section}
  
  \makeatother

\def\qed{{\hfill\hbox{\enspace${ \square}$}} \smallskip}
\def\sqr#1#2{{\vcenter{\vbox{\hrule height .#2pt \hbox{\vrule
 width .#2pt height#1pt \kern#1pt \vrule
width .#2pt} \hrule height .#2pt}}}}
\def\square{\mathchoice\sqr54\sqr54\sqr{4.1}3\sqr{3.5}3}

\def\ds{\begin{displaystyle}}
\def\eds{\end{displaystyle}}
\def\dis{\displaystyle }
\def\<{\langle }
\def\>{\rangle }

\def\R{\mathbb R}

\def\E{\mathbb E}
\def\P{\mathbb P}

\def\L{\mathbb L}
\def\D{\mathbb D}
\def\U{\mathbb U}

\def\cald{{\cal D}}

\def\calf{{\cal F}}
\def\calg{{\cal G}}
\def\calh{{\cal H}}

\def\calm{{\cal M}}

\def\calp{{\cal P}}
\def\calu{{\cal U}}

\def\call{{\cal L}}

\title{A Stochastic Optimal Control Problem
for the Heat Equation on the Halfline with Dirichlet Boundary-noise and Boundary-control }

\date{}
\author{
\\ \\
Federica Masiero\\
Dipartimento di Matematica e Applicazioni, Universit\`a di Milano-Bicocca\\
via R. Cozzi 53 - Edificio U5, 20125 MILANO, Italy\\
e-mail: federica.masiero@unimib.it}

\begin{document}

\maketitle
\begin{abstract}
 We consider a controlled state
equation of parabolic type on the halfline $(0,+\infty)$ with
boundary conditions of Dirichlet type in which the unknown 
is equal to the sum of the control and of a white
noise in time. We study finite horizon and infinite horizon 
optimal control problem related by menas of
backward stochastic differential equations.
\end{abstract}

\section{Introduction} \label{Introduction}
In this paper we study an optimal control problem for a state
equation of parabolic type on the halfline $(0,+\infty)$. We stress the fact that we
consider boundary conditions of Dirichlet type in which the unknown 
is equal to the sum of the control and of a white
noise in time, namely:
\begin{equation}\label{eqconcreta}
  \left\{
  \begin{array}{l}
  \dis
\frac{ \partial y}{\partial s}(s,\xi)= \frac{ \partial^2
y}{\partial \xi^2}(s,\xi)+f(s,y(s,\xi)), \qquad s\in [t,T],\;
\xi\in (0,+\infty),
\\\dis
y(t,\xi)=x(\xi),
\\\dis
 y(s,0)= u_s+\dot{W}_s, \quad
\end{array}
\right.
\end{equation}
In this equation $\{W_t,\;t\ge 0\}$ is a
standard real Wiener process,
$y(s,\xi,\omega)$ is the unknown real-valued process and represents 
the state of the system; the control is given by the 
real-valued processes $u(s,\xi,\omega)$ acting at $0$; $x:(0,+\infty)\rightarrow\R$.

Boundary control problems have been
widely studied in the deterministic literature (\cite{LasTri}) and have been addressed 
in the stochastic case as well (see\cite{DuMasl}, \cite{GoRoSw}, \cite{Ich}, \cite{Masl}).
In these works, the equation always contains noise also as 
a forcing term. 
In \cite{DebFuTe} a finite horizon optimal control problem for the 
stochastic heat equation with Neumann boundary conditions is
treated by backward stochastic differential equations. Here
we follow a similar approach but we consider the case with Dirichlet boundary conditions,
and we address both the finite horizon and 
the infinite horizon stochastic optimal control problems.
The main difficulties that we encounter in studying the control problem
for the state equation with Dirichlet boundary conditions are
related to the fact that the solution of equation (\ref{eqconcreta}) is not $L^2$-valued
unlike to the case of Neumann boundary conditions. Indeed, in \cite{DaZa-paper}
it is shown that, if we replace Neumann by
Dirichlet boundary conditions, the solution of (\ref{eqconcreta}) is
well defined in a negative Sobolev space $H^{\alpha}$, for $\alpha<-\dis\frac{1}{4}$.
Then in \cite{AlBo1}, see also \cite{BoGU}, it is shown that 
the solution $y(t,\cdot)$ of equation (\ref{eqconcreta}) with 
$u=0$ takes values in a weighted space $L^2((0,+\infty);\xi^{1+\theta}d \xi)$, 
nevertheless the problem was not reformulated as a stochastic evolution equation
in $L^2((0,+\infty);\xi^{1+\theta}d \xi)$.
The solutions are singular at the boundary, the singularity
is described in \cite{AlBo1} and \cite{Sow}. The reason is that the smoothing properties
of the heat equation are not strong enough to regularize a rough term such as a white noise.

In \cite{FaGo} equation (\ref{eqconcreta}), with $f=0$ is reformulated as 
an evolution equation in $L^2((0,+\infty);\xi^{1+\theta}d \xi)$
using results in \cite{Kr99} and in \cite{Kr01}. In these two 
papers it is shown that the 
Dirichlet Laplacian extends to a generator $A$ of 
an analytic semigroup on $L^2((0,+\infty);\xi^{1+\theta}d \xi)$.

%

Here we follow \cite{FaGo} and in Section 2
we reformulate equation (\ref{eqconcreta}) as a stochastic evolution
equation in $L^2((0,+\infty);\xi^{1+\theta}d \xi)$.
Namely we rewrite it as:
\begin{equation}\label{eqstatorifintro}
   \left\{
  \begin{array}{l}
  \dis
d X_s= AX_sds+F(s,X_s)ds+ BdW_s+Bu_sds, \qquad
s\in [t,T],
\\\dis
X_t=x,
\end{array}
\right.
\end{equation}
where $A$ stands for the Laplace operator with
homogeneous Dirichlet boundary conditions,
which is the generator of an analytic semigroup in $L^2((0,+\infty);\xi^{1+\theta}d \xi)$ 
(see \cite{Kr99} and \cite{Kr01}), $F$ is the evaluation operator
corresponding to $f$, $B=(\lambda-A)D_{\lambda}$ where $\lambda$ is an arbitrary positive
number and $D_{\lambda}$ is the Dirirchlet map 
(for more details on the abstract formulation of equation (\ref{eqconcreta}) see section 2.1).

The optimal control problem we wish to treat in this paper
consists in minimizing the following finite horizon cost
 \begin{equation}\label{costoconcreto}
J(t,x,u)=\E \int_t^T\int_0^{+\infty} \ell(s,\xi,
y(s,\xi),u_s)\;d\xi\;ds +\E \int_0^{+\infty} \phi(\xi,
y(T,\xi))\;d\xi.
\end{equation}

Our purpose is not only to prove existence of optimal controls but
mainly to characterize them by an optimal feedback law. To this aim first we
solve (in a suitable sense) the Hamilton-Jacobi-Bellman equation;
then we prove that such a solution is the value function
of the control problem and allows to
construct the optimal feedback law. Hamilton-Jacobi-Bellman equation
can be formally written as
\begin{equation}\label{HJBformale-intro}
  \left\{\begin{array}{l}\dis
\frac{\partial v(t,x)}{\partial t}+\call_t [v(t,\cdot)](x) = \Psi
(t, x,\nabla v(t,x)B),\quad t\in [0,T],\,
x\in H,\\
\dis v(T,x)=\Phi(x).
\end{array}\right.
\end{equation}
where $\call_t$ is the infinitesimal generator of the Markov
semigroup corresponding to the process $X$. We notice that 
$\call_t$ is highly degenerate, indeed
$\nabla^2f(x)$ appears only multiplied by $B$, and so
the equation \ref{HJBformale-intro} has very poor smoothing properties. 
%

We formulate the equation (\ref{HJBformale-intro}) in a mild sense, see for instance \cite{Go} and \cite{Go3}.
We notice that, when the state equation is linear, it is known
that the
semigroup $\{P_{s,t}[\,\cdot\,]: 0\leq s\leq t\}$ is strongly
Feller, nevertheless it seems that equation (\ref{HJBformale-intro})
cannot be solved by a fixed point argument as, for instance, in \cite{Go}
or \cite{Go3}, see also \cite{DebFuTe} and references therein.

We also mention here that, as it is well known, when the space is finite
dimensional Hamilton-Jacobi-Bellman equations can be successfully
treated using the notion of viscosity solution, see \cite{GoRoSw} for viscosity
approach to boundary optimal control. The point is
that, in the infinite dimensional case, very few uniqueness results
are available for viscosity solutions and all of them, obtained by
analytic techniques, impose strong assumptions on the operator $B$ and
on the nonlinearity $\Psi$, see, for instance, \cite{GoRoSw} or
\cite{Sw} and references within.

To solve the Hamilton-Jacobi-Bellman equation (\ref{HJBformale-intro}) in mild sense
we follow the approach based on Forward-Backward
stochastic differential equations, mainly developped, in a finite
dimensional setting,in the fundamental papers
\cite{kapequ}, \cite{PaPe-90} and \cite{PaPe}, and generalized, in infinite dimensions, 
in \cite{FuTe}.
The backward stochastic differential equation is in our case
\begin{equation}\label{BSDE-intro}
\left\{\begin{array}{l}
\dis dY^{t,x}_s= -\Psi(s,X_s^{t,x},Z^{t,x}_s)ds+Z^{t,x}_sdW_s,\quad s\in [t,T] \\
  Y^{t,x}_T=\Phi(X_T^{t,x})\\
\end{array}\right.
\end{equation}
and we need to study regular dependence of $Y$ on the 
initial datum $x$: in order to give sense to the term
$\nabla v(t,\cdot)B$ in \ref{HJBformale-intro} we have to 
differentiate $Y$ in the direction 
$(\lambda -A)^\alpha h$.

\noindent The control problem is solved by using the
probabilistic representation of the unique mild solution to
equation (\ref{HJBformale-intro}) which also gives existence of an
optimal feedback law, see Theorem \ref{th-rel-font}.

We also treat the infinite horizon optimal control problem:
minimize, over all admissible controls, the following infinite horizon cost
\begin{equation}\label{costoconcretoinfor}
J(x,u)=\E \int_0^{+\infty}e^{-\mu s}\int_0^{+\infty} \ell(s,\xi,
y(s,\xi),u_s)\;d\xi\;ds.
\end{equation}
The controlled state $y$ solves 
\begin{equation}\label{eqconcretaM}
  \left\{
  \begin{array}{l}
  \dis
\frac{ \partial y}{\partial s}(s,\xi)= \frac{ \partial^2
y}{\partial \xi^2}(s,\xi)-My(s,\xi)+f(s,y(s,\xi)), \qquad s\in [t,T],\;
\xi\in (0,+\infty),
\\\dis
y(t,\xi)=x(\xi),
\\\dis
 y(s,0)= u_s+\dot{W}_s, \quad
\end{array}
\right.
\end{equation}
where $M$ has to be taken sufficiently large,
see also lemma \ref{lemmader_forward_lim}.
The reason is that in the space $\calh$ we need to treat Dirichlet boundary conditions
it is not clear whether $A$ is dissipative or not.

%
As in the finite horizon case we consider mild solution of 
the Hamilton Jacobi Bellman equation related, which this time is
stationary.
The main tool for solving this stationary Hamilton-Jacobi-Bellman equation
is again BSDEs, where the final condition is replaced by boundedness
requirements on $Y$, see also \cite{FuTe-ell} and \cite{HuTess}..

\noindent As for the finite horizon case, in order to give sense to the term
$\nabla v(\cdot)B$, we have to differentiate
the following backward stochastic differential equation
\begin{equation}\label{BSDEintro_inf}
  dY^x_s=-\Psi(X_s^x,Z^x_s)\;ds+\mu Y^x_s\;ds+Z^x_s \;dW_s,\qquad s\geq 0,\\
\end{equation}
where $\mu>0$ and $\Psi$ is the hamiltonian function defined in a classical way.
To study the regularity property of equation (\ref{BSDEintro_inf}), we use similar ideas as in
\cite{HuTess}, where differentiability with respect to $x$ of $(Y^x,Z^x)$,
solution of an equation like (\ref{BSDEintro_inf}) with an arbitrary $\mu>0$, is investigated.
We notice again that since we have to give sense to $\nabla_xY^x(\lambda-A)^\alpha h$, for any
$h\in\calh$, we also need to differentiate equation (\ref{BSDEintro_inf}) in the
direction $(\lambda-A)^\alpha h$, and consequentely we have to study a
BSDE with some terms unbounded in time: such a situation is not studied in \cite{HuTess}.

The paper is structured as follows: in Section 2 we transpose the controlled
state equation in the infinite dimensional framework and we study
regularity properties of the solution of this (forward) state equation;
in Section 3 we study the backward equation
associated to the problem; in Section 4 we prove existence and
uniqueness of the Hamilton-Jacobi-Bellman partial differential
equation and in Section 5 we show how the previous results can be
applied to perform the synthesis of the optimal control, both in a strong 
and weak formulation.
Eventually we study the infinite horizon optimal control problem:
in Section 6 we study the regularity properties of the forward-bacward equations in infinite horizon,
in Section 7 we prove existence and uniquenes of the solution of
the stationary Hamilton-Jacobi-Bellman equation, and in section 8 we briefly present and solve the
infinite horizon optimal control problem.

\section{The forward equation}

In this section we introduce the ``concrete'' state equation, that we
reformulate in an abstract sense following \cite{FaGo}, and then we
study some regularity properties.

\subsection{Reformulation of the state equation}

We consider the following stochastic semilinear heat equation with control and noise on the boundary:
\begin{equation}\label{eqconcreta1}
  \left\{
  \begin{array}{l}
  \dis
\frac{ \partial y}{\partial s}(s,\xi)= \frac{ \partial^2
y}{\partial \xi^2}(s,\xi)+f(s,y(s,\xi)), \qquad s\in [t,T],\;
\xi\in (0,+\infty),
\\\dis
y(t,\xi)=x(\xi),
\\\dis
 y(s,0)= u_s+\dot{W}_s, \quad
\end{array}
\right.
\end{equation}
In this equation $\{W_t,\;t\ge 0\}$, is a
standard real Wiener process; 
$y(s,\xi,\omega)$ is the unknown real-valued process and represents 
the state of the system; the control is given by the 
real-valued process $u(s,\xi,\omega)$ which belongs to the class of 
admissible controls $\calu$, $f:[0,T]\times\R\to\R$ and $x:[0,+\infty)\to\R$.

It is our purpose to write the state equation as an evolution
equation in the space $\calh=L^2((0,+\infty);\xi^{1+\theta}d \xi)$, or in 
the space $L^2((0,+\infty);(\xi^{1+\theta}\wedge 1) d \xi)$, that we also denote by $\calh$. 
The parameter $\theta \in (0,1)$.
\noindent On equation (\ref{eqconcreta1}) we assume that 
\begin{hypothesis}\label{ipotesiconcrete}
\begin{enumerate}
\item[1)] The function $f:[0,T]\times\R\to\R$ is measurable,
for every $t\in [0,T]$ the function
  $f(t,\cdot):\mathbb{R} \rightarrow \mathbb{R}$
  is continuously differentiable and there exists a
constant $C_f$ such that
$$
|f(t,0)|+\left|\frac{\partial f}{\partial r}(t,r)\right|\le C_f,
\qquad t\in[0,T],\;r\in\R.
$$
\item[2)] The initial condition $x(\cdot)$ belongs to $\calh$.
\item[3)] The set of admissible control actions $\calu$
is a bounded closed subset of $\R$.
\end{enumerate}
\end{hypothesis}
Equation \ref{eqconcreta}, in the case of $f=0$, is reformulated as an evolution equation 
in $\calh$ in \cite{FaGo} and we follow that approach. 
Let us denote by $A$ the Laplacian operator with Dirichlet boundary conditions:
it is proved in \cite{Kr01} that the strongly continuous 
heat semigroup generated in $L^2((0,+\infty))$ by $A$
extends to a bounded $C_0$ semigroup $(e^{tA})_{t\geq 0}$ in $\calh$ with generator 
 still denoted by $A:\cald(A)\subset \calh \rightarrow \calh$. 
The semigroup $(e^{tA})_{t\geq 0}$ is analytic. So, for every $\beta>0$,
\begin{equation}
 \label{stima-analitica}
\Vert (\lambda-A)^\beta e^{tA} \Vert\leq C_\beta t^{-\beta} \qquad \text{for all} \qquad t\geq 0.
\end{equation}

Let us also introduce the Dirichlet map: for given $\lambda>0$, let 
$D_\lambda:\R\rightarrow \calh$ be such that, for $a\in \R$, $D_\lambda (a)=a\psi_\lambda$, where 
$\psi_\lambda:\R\rightarrow\R^+$, $\psi_\lambda: \xi \mapsto e^{-\lambda\xi}.$
In the following proposition we collect some results contained in \cite{FaGo}. From now on 
$\lambda>0$ is fixed.
\begin{proposition}\label{propFabbriGoldys}
For all $\alpha\in[0,\frac{1}{2}+\frac{\theta}{4})$, $\psi_\lambda\in D((\lambda-A)^\alpha)$ ,
and in particular $D_\lambda\in \call(\R;D((\lambda-A)^\alpha))$. So the operator
\begin{equation}\label{operatoreB}
 B:=(\lambda-A)D_\lambda:\R\rightarrow \calh^{\alpha-1}
\end{equation}
is bounded, and for every $t>0$ the operator 
\begin{equation}\label{operatoreB_bis}
 (\lambda-A)e^{tA}D_\lambda=(\lambda-A)^{1-\alpha}e^{tA}(\lambda-A)^\alpha D_\lambda:\R\rightarrow \calh
\end{equation}
is bounded as well. From now on let  
$\alpha\in(\frac{1}{2},\frac{1}{2}+\frac{\theta}{4})$.
For all $\gamma<2\alpha-1$, the following holds
\begin{itemize}
 \item [i)] For each $t>0$, the operator $e^{tA}B:\R\rightarrow\calh$
   is bounded and the function $t\mapsto e^{tA}Ba$ is continuous $\forall$ $a\in\R$.
 \item [ii)] \begin{equation}
\label{stimaHS}
\int_0^T s^{-\gamma}\Vert e^{sA}B \Vert_{HS}^2 ds <+\infty
\end{equation}
where $\Vert \cdot \Vert_{HS}$ stands for the Hilbert-Schmidt norm.
 \item [iii)] For every $0\leq t <T$ the stochastic convolution
 \begin{equation}
\label{convstoc}
W_A(s)=\int_t^s e^{(s-r)A}B dW_r, \qquad s\in[t,T]
\end{equation}
is well defined, belongs to $L^2(\Omega;C([t,T],\calh))$ and 
has continuous trajectories in $\calh$. 
 \item[iv)] For every $0\leq t <T$ 
and $u\in L_\calp^2(\Omega\times[0,T];\R)$
 \begin{equation}
\label{int_controllo}
I_s=\int_t^s e^{(s-r)A}B u_r dr, \qquad s\in[t,T]
\end{equation}
is well defined in $L_\calp^2(\Omega\times[0,T];\R)$. 
Moreover $I\in L^2(\Omega;C([t,T],\calh))$ and 
$\Vert I \Vert_{ L^2(\Omega;C([t,T],\calh))}\leq C 
\Vert u \Vert _{L_\calp^2(\Omega\times[0,T];\R)}$.
\end{itemize}
\end{proposition}
\noindent We want to rewrite equation (\ref{eqconcreta1}) as en evolution
equation in $\calh$. The state will be
denoted by $X^u_s=y(s,\cdot)$. Thus $\{X^u_s, s\in [t,T]\}$ is a
process in $\calh$ and the initial condition is assumed to belong to
$\calh$. Equation (\ref{eqconcreta1}), in the case of $f=0$, can now be reformulated as
\begin{equation}\label{eqstatoformalelineare}
  \left\{
  \begin{array}{l}
 dX^u_s= AX^u_s ds+Bu_s ds +B dW_s  \qquad s\in [t,T],
\\\dis
X^u_t=x,
\end{array}
\right.
\end{equation}

\begin{definition}\label{defsolmildlineare}
An $\calh$-valued predictable process $X$  is called a mild solution 
to equation (\ref{eqstatoformalelineare}) on $[0,T]$ if 
\[
 \P\int_0^T\vert X^u_r\vert^2 dr <+\infty
\]
and, for every $0\leq t<T$, $X$ satisfies the integral equation
\[
 X^u_s=e^{(s-t)A}+\int_t^s e^{(s-r)A} Bu_r dr +\int_t^s e^{(s-r)A} BdW_r 
\]
\end{definition}
Following \cite{FaGo}, theorem 2.6, we state the following:
\begin{theorem}\label{teolineare}
Assume that hypothesis \ref{ipotesiconcrete} holds true, 
then equation (\ref{eqstatoformalelineare}) has, according 
to definition \ref{defsolmildlineare}, a 
unique mild solution $X\in L^2(\Omega;C([t,T],\calh))$. Moreover if 
$u=0$ then $X$ is a Markov process in $\calh$.
\end{theorem}

Next we want to give an abstract reformulation in $\calh$ 
of the semilinear equation (\ref{eqconcreta1}). 
We define $F:[0,T]\times \calh\rightarrow \calh$ setting for $s\in [0,T]$
and $X\in\calh$
\begin{equation}\label{notazioniF}
F(s,X)(\xi)=f(s,X(\xi)) 
\end{equation}
By hypothesis \ref{ipotesiconcrete}, point 1), it turns out that
$F:[0,T]\times \calh\to \calh$ is a measurable
function and
$$
|F(t,0)| + |F(t,x_1)-F(t,x_2)|\le C_f (1+|x_1-x_2|),
\qquad t\in [0,T],\;x_1,x_2\in \calh.
$$
Moreover, for every $t\in [0,T]$, $F(t,\cdot)$ has a G\^{a}teaux derivative
$\nabla_x F(t,x)$ at every point $x\in \calh$, and we get that $|\nabla F(t,x)|\le C_f$.
Finally, the function $(x,h)\to \nabla F(t,x)h$ is continuous
as a map $\calh\times\calh \to \R$.
Note that we consider $\nabla F(t,x)$ as an element of $\calh^*$ and
we denote its action on $h\in \calh$ by $\nabla F(t,x)h$.
Equation (\ref{eqconcreta1}) can now be reformulated as
\begin{equation}\label{eqstatoformale}
  \left\{
  \begin{array}{l}
  \dis
dX^u_s= AX^u_sds+F(s,X^u_s)ds+Bu_s ds +B dW_s  \qquad s\in [t,T],
\\\dis
X^u_t=x,
\end{array}
\right.
\end{equation}
The equation (\ref{eqstatoformale}) is formal. The precise meaning
of the state equation is in the following
\begin{definition}\label{defsolmild}
An $\calh$-valued predictable process $X$  is called a mild solution 
to equation (\ref{eqstatoformale}) on $[0,T]$ if 
\[
 \P\int_0^T\vert X^u_r\vert^2 dr <+\infty
\]
and for every $0\leq t<T$, $X$ satisfies the integral equation
\begin{equation}
\label{eqstatocontrol}
 X^u_s=e^{(s-t)A}x+\int_t^s e^{(s-r)A}F(r,X^u_r)\;dr 
+\int_t^s e^{(s-r)A} Bu_r dr +\int_t^s e^{(s-r)A} BdW_r.
\end{equation}
\end{definition}
We now prove existence and uniqueness of a mild solution of equation (\ref{eqstatoformale})
\begin{theorem}\label{teosolmild}
Assume that hypothesis \ref{ipotesiconcrete} holds true, then equation 
(\ref{eqstatoformale}) has, according to definition \ref{defsolmild}, a 
unique mild solution $X\in L^2(\Omega;C([t,T],\calh))$ and, if 
$u=0$, $X$ is a Markov process in $\calh$.
\noindent Moreover for every $p\in
[1,\infty)$, $\alpha\in [0,\theta/4)$, $t\in [0,T]$ there exists a
constant $c_{p,\alpha}$ such that
\begin{equation}\label{stimaes}
    \E\sup_{s\in (t,T]}(s-t)^{p\alpha}|X_s^{t,x}|^p_{D(-A)^\alpha}\leq
c_{p,\alpha}(1+|x|_\calh)^p.
\end{equation}

\end{theorem}
{\bf Proof.} 
We consider the Picard approximation scheme; for the sake of
simplicity we consider $u=0$ in equation (\ref{eqstatoformale}) 
and we denote by $(X_s)_{s\in[t,T]}$ the solution. We define
\begin{align*}
& X_s^0=e^{(s-t)A}x,\\ \nonumber
& X_s^{n+1}=e^{(s-t)A}x+\int_t^s e^{(s-r)A}F(r,X^n_r)\;dr 
 +\int_t^s e^{(s-r)A} BdW_r, \qquad n\geq 0.\\ \nonumber
\end{align*}
By induction it follows that for every $n\geq 0$, $X^n \in L^2(\Omega;C([t,T],\calh)$.
Moreover, by equipping $L^2(\Omega;C([t,T],\calh)$ with the equivalent norm 
\[
 \Vert Y \Vert^2_{\beta,L^2(\Omega;C([t,T],\calh)}=\E \sup_{s\in[t,T]}e^{-\beta s}\vert Y_s \vert^2_\calh
\]
it turns out that $(X^n)_n$ is a Cauchy sequence in
$(L^2(\Omega, C([t,T],\calh)), \Vert \cdot \Vert^2_{\beta,L^2(\Omega;C([t,T],\calh)})$, 
whose limit is the unique mild solution to equation (\ref{eqstatoformale}).
Next we want to prove estimate (\ref{stimaes}): first we prove that the 
stochastic convolution defined in (\ref{convstoc}) belongs to 
$L^2(\Omega;C([t,T],D(-A)^\alpha))$.
By the factorization method, see e.g. \cite{DP1}, p. 128, let $\gamma \in (0,\frac{1}{2})$:
the stochastic convolution can be written as 
\begin{equation*}
 W_A(s)=\frac{\sin \pi \gamma}{\pi} \int_t^s e^{(s-r)A}(s-r)^{\gamma-1}Y_rdr,
\end{equation*}
where 
\begin{equation*}
 Y_r= \int_t^r (r-\sigma)^{-\gamma}e^{(r-\sigma)A}BdW_r.
\end{equation*}
Let us write, for $\beta\in(\frac{1}{2},\frac{1}{2}+\frac{\theta}{4} )$, 
$e^{(r-\sigma)A}B=(\lambda-A)^{1-\beta}e^{(r-\sigma)A}(\lambda-A)^\beta D_\lambda$,
so that the stochastic convolution is given by
\begin{equation*}
 W_A(s)=\frac{\sin \pi \gamma}{\pi} \int_t^s e^{(s-r)A}(s-r)^{\gamma-1}
(\lambda-A)^{1-\beta}\widehat{Y}_rdr,
\end{equation*}
where 
\begin{equation*}
 \widehat{Y}_r= \int_t^r (r-\sigma)^{-\gamma}(\lambda-A)^\beta D_\lambda dW_r.
\end{equation*}
It turns out that $\widehat{Y}\in L^p(\Omega\times [t,T],\calh)$, and so, see 
e.g. \cite{DP2} Proposition A.1.1, for $\alpha+1-\beta+\frac{1}{p}<\gamma<\frac{1}{2}$, 
$(\lambda-A)^\alpha W_A\in L^p(\Omega, C ([t,T],\calh) $, where
\begin{equation*}
 (\lambda-A)^{\alpha}W_A(s)=\frac{\sin \pi \gamma}{\pi} \int_t^s e^{(s-r)A}(s-r)^{\gamma-1}
(\lambda-A)^{\alpha +1-\beta}\widehat{Y}_rdr.
\end{equation*}
We can conclude that 
for $\alpha+1-\beta+\frac{1}{p}<\gamma<\frac{1}{2}$, i.e. for 
$\alpha<\frac{\theta}{4}$, and $p$ sufficiently large, 
$ W_A\in L^p(\Omega, C ([t,T],D(\lambda-A)^\alpha) $.
In a similar, and simpler way, if $u\neq 0$,
we could treat the term 
\begin{equation*}
  \int_t^s e^{(s-r)A}B u_r dr.
\end{equation*}

For $a>0$ we denote by $\mathbb{K}_{a,\alpha,t}$ the Banach space of
all predictable processes $X: \Omega \times (t,T]\rightarrow
D(\lambda-A)^\alpha$
such that
$$|X|^p_{\mathbb{K}_{a,\alpha,t}}:= \E \sup_{s\in (t,T]}
e^{pa s}(s-t)^{p\alpha}|X_{s}|^p_{D(\lambda-A)^\alpha}< +\infty$$
endowed with the
above norm. We have just shown that $W_A \in
\mathbb{K}_{a,\alpha,t}$. Moreover, for all $x\in H$,
$$\sup_{s\in (t,T]} (s-t)^{\alpha} |e^{(s-t)A}x|_{D(\lambda-A)^\alpha}\leq c|x|.$$
Thus if we define for $X\in \mathbb{K}_{a,\alpha,t}$
$$\Lambda(X,t)(s)=\int_t^s
e^{(s-r)A}F(r,X_r)\;dr+e^{(s-t)A}x+W_A(s),$$ it is
immediate to prove that $\Lambda(X,t)\in \mathbb{K}_{a,\alpha,t}$.
Moreover by straightforward estimates
$$ |\Lambda(X^1,t)- \Lambda(X^2,t)|^p_{\mathbb{K}_{a,\alpha,t}}\leq
g^p(a)C_ F^p |X^1-X^2|^p_{\mathbb{K}_{a,\alpha,t}}$$ where
$$g(a)=\sup_{t\in [0,T]}t^{1-\alpha}\int_0^1
(1-s)^{-\alpha}s^{-\alpha}e^{-ats}ds.
$$
By the Cauchy-Schwartz inequality $g(a)\leq
T^{1/2-\alpha}a^{-1/2}\left(\int_0^1
(1-s)^{-2\alpha}s^{-2\alpha}ds\right)^{1/2}$
 thus if $a$ is large enough $\Lambda(\cdot,t)$ is a
contraction in $\mathbb{K}_{a,\alpha,t}$. The unique fixed point
is clearly a mild solution of equation (\ref{eqstatoformale}) and
(\ref{stimaes}) holds. Uniqueness is an immediate consequence of
the Gronwall lemma.
\qed

It is also useful to consider
the uncontrolled version of equation
(\ref{eqstatoformale}) namely:
\begin{equation}\label{eqstatonocontrol}
X_s=\dis e^{(s-t)A}x+\int_t^s e^{(s-r)A}F(r,X_r)\;dr +
\int_t^se^{(s-r)A}B\;dW_r, \qquad s\in [t,T].
\end{equation}
We will refer to (\ref{eqstatonocontrol}) as the forward equation.

\subsection{Regular dependence on initial conditions.}

In this section we consider again the solution
of the forward equation (\ref{eqstatonocontrol}),
i.e. of the uncontrolled state
equation on the time interval $[t,T]$ with initial condition
$x\in \calh$. It will be denoted by $X_s^{t,x}$,
 to stress dependence on the initial
data $t$ and $x$.
It is also convenient
to extend the process $X_\cdot^{t,x}$ letting $X_s^{t,x}=x$ for
$s\in [0,t]$. In a similar way, we extend also the stochastic
convolution by setting $W_A(s)=0$ for $s\in [0,t)$.
From now on we assume that Hypothesis \ref{ipotesiconcrete}
holds.

\noindent We study the dependence of the process
$\{X_s^{t,x},\;s\in [0,T]\}$ on the parameters $t,x$.

\begin{proposition}\label{regequazforward}
For
any $p\geq 1$ the following holds:
\begin{enumerate}
    \item the map $(t,x)\to X_\cdot^{t,x}$
defined on $[0,T]\times \calh$ and with
  values in $L^p_\calp (\Omega;C([0,T];\calh))$ is continuous.
    \item For every $t\in [0,T]$ the map $x\to X_\cdot^{t,x}$
has, at every point $x\in \calh$, a G\^{a}teaux derivative $\nabla_x
X_\cdot^{t,x}$. The map $(t,x,h)\to \nabla_x X_\cdot^{t,x}h$ is
continuous as a map $[0,T]\times \calh\times \calh\to L^p_\calp
(\Omega;C([0,T];\calh))$ and, for every $h\in \calh$, the following
equation holds $\P$-a.s.:
\begin{equation}\label{derdixmild}
    \nabla_x X_s^{t,x}h=e^{(s-t)A}h+
    \int_t^s e^{(s-\sigma)A}\nabla_xF(\sigma,
    X_\sigma^{t,x}) \nabla_x X_\sigma^{t,x}\; d\sigma,
    \qquad s\in [t,T],
\end{equation}
and $\nabla_x X_s^{t,x}h =h$ for $s\in [0,t]$.

\end{enumerate}
\end{proposition}

{\bf Proof.} 
We start by proving continuity. 
We begin considering the stochastic convolution:
we know that $\int_t^se^{(s-r)A}BdW_r=W_A(s)\in
L^p_{\mathcal{P}}(\Omega,C([0,T];\calh))$ and we have to prove that the map
$t\rightarrow \int_t^se^{(s-r)A}BdW_r$ is continuous with values in
$L^p_{\mathcal{P}}(\Omega,C([0,T];\calh))$. Fix $t\in [0,T]$, 
$\beta \in (\frac{1}{2}, \frac{1}{2}+\frac{\theta}{4})$ and let
$t_n\rightarrow t^+$, (in a similar way if $t_n\rightarrow t^-$)
\begin{align*}
\E\sup_{s\in [t,T]}\left\vert\int_t^{s\wedge
t_n} e^{(s-\sigma)} B dW_{\sigma}\right\vert_\calh^p &
\leq   \E\sup_{s\in [t,T]}\left|\int_t^{s\wedge
t_n}(\lambda-A)^{1-\beta}e^{(s-\sigma)A}
(\lambda-A)^{\beta}D_\lambda dW_{\sigma}\right|_\calh ^p\\ \nonumber
&\leq \sup_{s\in [t,T]}  \left(\int_t^{s\wedge
t_n}\left|(\lambda-A)^{1-\beta}e^{(s-\sigma)}
(\lambda-A)^{\beta}D_\lambda\right|^2_\calh d\sigma\right)^\frac{p}{2}\\ \nonumber
& \leq C\sup_{s\in [t,T]}\left(\int_t^{s\wedge
t_n}(s-\sigma)^{2(1-\beta)}d\sigma\right)^{p/2}\rightarrow 0.
\end{align*}
Similarly if we extend $e^{(s-t)A}x=x$ for $s<t$ then
$$
\sup_{s\in [0,T]}\left|e^{(s-t_n)A}x-e^{(s-t)A}x\right|
\rightarrow 0
$$
as $t_n\rightarrow t$; moreover the map $x\rightarrow e^{(\cdot-t)A}x$ considered with
values in $C([0,T],\calh)$ is clearly continuous in $x$ uniformly in
$t$.

\noindent Now let $t_n\rightarrow t^+$ and $x_n\rightarrow x$: 
\begin{align*}
& \E\sup_{s\in [0,T]}|X_s^{t_n,x_n}-X_s^{t,x}|_\calh^p \\ \nonumber
&\leq C\sup_{s\in [0,T]}\left|e^{(s-t_n)A}x-e^{(s-t)A}x\right|_\calh^p 
+ C\E\sup_{s\in [0,T]}\left\vert\int_t^{s\wedge
t_n} e^{(s-\sigma)} B dW_{\sigma}\right\vert_\calh^p \\ \nonumber
&+ C\E\sup_{s\in [0,T]} \left\vert\int_{t_n}^s e^{(s-\sigma)A}F(\sigma,X_\sigma^{t_n,x_n})d\sigma
-\int_t^s e^{(s-\sigma)A}F(\sigma,X_\sigma^{t,x})d\sigma\right\vert^p_\calh \\ \nonumber
& \leq \epsilon (\vert x_n-x\vert_\calh, \vert t_n-t \vert)+
+C_{F,T}\int_{t_n}^t \vert X_\sigma^{t,x} - X_\sigma^{t_n,x_n}\vert^p d\sigma 
+C_T (t_n-t)^{1/p}(1+\vert x\vert_\calh^p), \\ \nonumber
\end{align*}
where $\epsilon (\vert x_n-x\vert_\calh, \vert t_n-t \vert):=
\sup_{s\in [0,T]}\left|e^{(s-t_n)A}x-e^{(s-t)A}x\right|+
\E\sup_{s\in [0,T]}\left\vert\int_t^{s\wedge
t_n} e^{(s-\sigma)} B dW_{\sigma}\right\vert_\calh^p $
and $C_{F,T}$ is a constant that depends on $F$ and on $T$.
By the Gronwall lemma
\[
 \E(\sup_{s\in [0,T]}|X_s^{t_n,x_n}-X_s^{t,x}|_\calh^p \rightarrow 0
\]
as $t_n\rightarrow t^+$ and $x_n\rightarrow x$.

\noindent The proof of differentiability is similar to the proof of proposition 3.1 in
\cite{DebFuTe}, and we omit it.
\qed

\begin{proposition}\label{regequazforward2} For every $\alpha \in
[0,1)$ there exists a family of predictable processes
$\{\Theta^{\alpha}(\cdot,t,x)h: h\in H,\,x\in H,\, t\in [0,T]\} $
all defined on $\Omega\times[0,T]\to H$ such that the following
holds:
\begin{enumerate}
\item the map $h \rightarrow \Theta^{\alpha}(\cdot,t,x)h$ is
linear and, if $h\in D(\lambda -A)^\alpha$, then
\begin{equation}\label{defdiTheta}
    \Theta^{\alpha}(s,t,x)h=\left\{\begin{array}{ll}
       \dis \left(\nabla_x X_s^{t,x}- e^{(s-t)A}\right)(\lambda -A)^\alpha h
       &\hbox{ if } s\in [t,T], \\
      0 &\hbox{ if } s\in [0,t).
    \end{array}\right.
\end{equation}
    \item the map $(t,x,h)\to\Theta^{\alpha}(\cdot,t,x)h$ is
continuous $[0,T]\times \calh\times \calh\to L^\infty _\calp
(\Omega;C([0,T];\calh))$.
    \item there exists a constant $C_{\theta,\alpha}$ such that
\begin{equation}\label{stimadi nabla theta}
|{\Theta}^{\alpha}(\cdot,t,x)h|_{L^{\infty}_{\mathcal{P}}(\Omega,C([0,T];\calh))}\leq
C_{\theta,\alpha} |h| \hbox{ for all $t\in [0,T]$, $x,h\in \calh$}.
\end{equation}
\end{enumerate}
\end{proposition}

{\bf Proof:} For fixed $t\in [0,T]$ and $x,h\in \calh$ consider the
equation:
\begin{equation}\label{equazpertheta}
    \begin{array}{rcl}
   \dis {\Theta}^{\alpha}(s,t,x)h&=& \dis \int_t^s e^{(s-\sigma)A}\nabla_x
F(\sigma,X_\sigma^{t,x}){\Theta}^{\alpha}(\sigma,t,x)h
d\sigma\\
\quad \quad \dis && \dis +\int_t^s e^{(s-\sigma)A} \nabla_x
F(\sigma,X_\sigma^{t,x})(\lambda-A)^{\alpha}e^{(\sigma-t)A}h\;
d\sigma.
\end{array}
\end{equation}
Notice that $$\int_t^s\left|e^{(s-\sigma)A} \nabla_x
F(\sigma,X_\sigma^{t,x})(\lambda-A)^{\alpha}e^{(\sigma-t)A}h
\right|d\sigma\leq C_f \int_t^s(\sigma-t)^{-\alpha}|h|d\sigma\leq
c|h|$$ for a suitable constant $c$.

Since $\nabla_x F$ bounded it is immediate to prove
that equation (\ref{equazpertheta}) has $\mathbb{P}$-almost surely
a unique solution in $C([t,T];\calh)$. Moreover extending
${\Theta}^{\alpha}(s,t,x)h=0$ for $s<t$ and considering it as a
process we have ${\Theta}(\cdot,t,x)h \in
L^{\infty}_{\mathcal{P}}(\Omega,C([0,T];\calh))$ and
$|{\Theta}^{\alpha}(\cdot,t,x)h|_{L^{\infty}_{\mathcal{P}}(\Omega,C([0,T];\calh))}\leq
C_{\alpha} |h|$. The continuity with respect to $t$, $x$ and
$h$ can be easily shown as in the proof of the previous Proposition. Moreover
linearity in $h$ is straight-forward. Finally for all $k\in
D(\lambda-A)^{\alpha}$ setting $h=(\lambda-A)^{\alpha}k$ equation
(\ref{derdixmild}) can be rewritten: $$
\begin{array}{l}
   \dis \left(\nabla_x
   X(s,t,x)(\lambda-A)^{\alpha}k-e^{(s-t)A}(\lambda-A)^{\alpha}k\right)=\\
   \quad \quad \dis
   +\int_t^s e^{(s-\sigma)A} \nabla_x
F(\sigma,X(\sigma,t,x))e^{(\sigma-t)A}(\lambda-A)^{\alpha}k
d\sigma \\
\quad \quad \dis +\int_t^s e^{(s-\sigma)A} \nabla_x
F(\sigma,X(\sigma,t,x))\left( \nabla_x
X(\sigma,t,x)(\lambda-A)^{\alpha}k-e^{(\sigma-t)A}(\lambda-A)^{\alpha}k\right)
d\sigma.
\end{array} $$
Comparing the above equation with equation (\ref{equazpertheta})
by the Gronwall Lemma we get
 $\Theta^{\alpha}(s,t,x)k=
       \dis \left(\nabla_x X(s,t,x)- e^{(s-t)A}\right)(\lambda -A)^\alpha
       k$ $\mathbb{P}$-a.s. for all
       $ s\in [t,T]$.
 \qed

\subsection{Regularity in the Malliavin sense.}

In order to state the following results we need to
recall some basic definitions from the Malliavin calculus,
mainly to fix notation.
We refer the reader to the book
\cite{Nu} for a detailed exposition; the paper
\cite{GrPa} treats the extensions to
Hilbert space valued random variables and processes.

For every $h\in L^2([0,T];\R)$ we denote
$$
W(h)=\int_0^Th(t)\; dW(t).
$$
Given a Hilbert space $K$,
let $S_K$ be the set of $K$-valued random variables $F$ of the
form
$$
F=\sum_{j=1}^mf_j(W(h_1),\ldots,W(h_n))e_j,
$$
where
 $h_1,\ldots , h_n\in L^2([0,T];\R)$,
  $\{e_j\}$ is a basis of $K$ and $f_1,\ldots f_m$ are
infinitely differentiable functions $\R^n\to \R$ bounded together
with all their derivatives.
The Malliavin derivative $DF$ of $F\in S_K$ is defined as the
process $\left\lbrace D_s, s\in [0,T]\right\rbrace $, where
$$D_sF= \sum_{j=1}^m\sum_{k=1}^n\partial_k f_j
(W(h_1),\ldots,W(h_n))\; h_k^i(s)\;e_j.
$$
By $\partial_k$ we denote the partial
derivative with respect to the $k$-th variable.
$DF$ is a process
with values in $K$, that we will identify with
an element of $L^2(\Omega\times [0,T]; K)$
with the norm:
$$
\|DF\|_{L^2(\Omega\times [0,T]; K)}^2
= \E\int_0^T \|D_sF\|^2_{K}\;ds.
$$
 It is known that
the operator
$D:S_K\subset L^2(\Omega;K)\to L^2(\Omega\times [0,T]; K)$
is closable. We denote by $\D^{1,2}(K)$ the domain
of its closure, and use the same letter to denote
$D$ and its closure:
$$
D:\D^{1,2}(K) \subset L^2(\Omega;K)
\to L^2(\Omega\times [0,T]; K).
$$
The adjoint operator of $D$,
$$
\delta: {\;\rm dom\;}(\delta)\subset
L^2(\Omega\times [0,T]; K)\to L^2(\Omega;K),
$$
is called Skorohod integral. It is known that
dom$(\delta)$ contains
$L_\calp^2(\Omega;L^2([0,T]; K))$ and the Skorohod
integral of a process in this space
coincides with the It\^o integral. The class $\L^{1,2}(K)$
is also contained in dom$(\delta)$ ,
the latter being defined as the space of processes
$u\in L^2(\Omega\times [0,T]; K)$ such that
$u_r\in \D^{1,2}(K)$ for a.e. $r\in [0,T]$
 and there exists a
measurable version of $D_su_r$ satisfying
$$
\| u \|^2_{\L^{1,2}(K)}=
 \| u \|^2_{L^2(\Omega\times [0,T]; K)}
+\E\sum_{i=1}^2\int_0^T \int_0^T \|D^i_su_r\|^2_{K}\,dr\,ds<\infty.
$$
Moreover, $\| \delta(u) \|^2_{ L^2(\Omega;K)}\leq
\| u \|^2_{\L^{1,2}(K)}$.
The definition of $\L^{1,2}(K)$ for an arbitrary Hilbert space $K$
is entirely analogous; clearly,
$\L^{1,2}(K)$ is isomorphic to $L^2([0,T]; \D^{1,2}(K))$.

Finally we recall that if $F\in \D^{1,2}(K)$ is $\calf_t$-adapted
then $DF=0$ a.s. on $\Omega
\times (t,T] $.

Now for $(t,x)$ fixed let us consider again the process
$\{ X_s^{t,x}, s\in [t,T]\}$ solution of the
forward equation
(\ref{eqstatonocontrol}). It will be
denoted simply by $\{ X_s, s\in [t,T]\}$ or even $X$.
 We still agree that $X_s=x$ for $s\in [0,t)$.
We will soon prove that $X$
  belongs to $\L^{1,2}(\calh)$. Then it is clear that the equality
  $D_\sigma X_s=0$ $\P$-a.s. holds
  for a.a. $\sigma,t,s$ if $s<t$ or
$\sigma>s$.

In the rest of this section we still assume that
Hypothesis \ref{ipotesiconcrete}
holds.

\begin{proposition}\label{regequazforwardmalliavin}
Let $t\in [0,T]$ and $x\in \calh$ be fixed. Then
$X\in
 \L^{1,2}(\calh)$, and $\P$-a.s. we have, for a.a.
$\sigma, s$ such that $t\le \sigma\le s\le T$, and for 
$\beta\in (0,\frac{1}{2}+\frac{\theta}{4})$
\begin{equation}\label{derivmall}
    D_\sigma X_s= (\lambda - A)\,e^{(s-\sigma)A}b+\int_\sigma^s
e^{(s-r)A}\nabla F(r,X_r)\;D_\sigma X_r\;dr,
\end{equation}
\begin{equation}\label{stimasuderivmall}
  | D_\sigma X_s|
  \le C(s-\sigma)^{\beta-1}.
\end{equation}
Moreover for every $s\in [0,T]$ we have
$X_s\in \D^{1,2}(\calh)$ and
$DX_s\in L^\infty (\Omega;L^2([0,T]; \calh))$.

Finally, for every $q\in [2,\infty)$
the map $s\to X_s$ is continuous from $[0,T]$ to $L^q(\Omega;\calh)$
and the map $s\to DX_s$ is continuous from $[0,T]$ to
$L^q(\Omega;L^2([0,T]; \calh))$.
\end{proposition}

{\bf Proof.}
For simplicity of notation we write the proof for the case
$t=0$. Thus,
\begin{equation}\label{eqdibase}
X_s=\dis e^{sA}x+\int_0^s e^{(s-r)A}F(r,X_r)\;dr +
\int_0^se^{(s-r)A}B\;dW_r \qquad s\in [0,T].
\end{equation}
We set $J_n=n(n-A)^{-1}$ and we consider the approximating equation
\begin{equation}\label{eqdibaseapprox}
\begin{array}{l}\dis
 X^n_s=\dis e^{sA}x+\int_0^s e^{(s-r)A}F(r,X^n_r)\;dr +
\int_0^se^{(s-r)A}J_nB\;dW_r,
\\\dis
\qquad \qquad
=\dis e^{sA}x+\int_0^s e^{(s-r)A}F(r,X^n_r)\;dr +
\int_0^se^{(s-r)A}(\lambda
-A)J_nD_\lambda\;dW_r,
\qquad s\in [0,T].
\end{array}
\end{equation}
Since $(\lambda
-A)J_n$ is a linear bounded operator in $\calh$, we can
apply Proposition 3.5 of \cite{FuTe}
and conclude that $X^n\in \L^{1,2}(\calh)$, and
that $\P$-a.s. we have, for a.a.
$\sigma, s$ such that $0\le \sigma\le s\le T$,
\begin{equation}\label{derivmallappr}
    D_\sigma X^n_s=
  (\lambda - A)\,e^{(s-\sigma)A}J_nD_\lambda+\int_\sigma^s
e^{(s-r)A}\nabla F(r,X^n_r)\;D_\sigma X^n_r\;dr.
\end{equation}
Since for $0<\beta<\frac{1}{2}+\frac{\theta}{4}$
$$
| (\lambda - A)\,e^{(s-\sigma)A}J_nD_\lambda|
\le
|(\lambda -A)^{1-\beta}e^{(s-r)A}| \,
|J_n|\, |(\lambda-A)^{\beta}D_\lambda|
\le
C(s-r)^{\beta-1},
$$
by the boundedness of $\nabla F$ and the Gronwall lemma
it is easy to deduce that
$|D_\sigma X^n_s|\le C(s-\sigma)^{\beta-1}$. In particular
it follows that $DX^n$ is bounded
in the space $L^2(\Omega\times [0,T]\times [0,T];\calh)$.

Subtracting (\ref{eqdibaseapprox}) from
(\ref{eqdibase}) and using the Lipschitz
character of $F$ we obtain
$$
\E| X^n_s-X_s|^2\le C\int_0^s
\E| X^n_r-X_r|^2\;dr +C
\int_0^s|(\lambda -A)^{1-\beta}e^{(s-r)A}(\lambda
-A)^{\beta}(J_nD_\lambda-D_\lambda)|^2\;dr.
$$
For $\beta>\frac{1}{2}$ the last integral can be estimated by
$$
C
\int_0^s(s-r)^{2\beta-2}\;dr
|(J_n-I)(\lambda-A)^{\beta}D_\lambda|^2
\le C|(J_n-I)(\lambda-A)^{\beta}D_\lambda|^2,
$$
and so as $n\rightarrow \infty $ it tends to zero for , by well-known properties of the operators
$J_n$. If follows from the Gronwall lemma that
$\sup_s\E| X^n_s-X_s|^2\to 0$ and in particular
$X^n\to X$ in $L^2(\Omega\times [0,T];\calh)$.

The boundedness of the sequence $DX^n$ proved before
and the closedness of the operator $D$ imply that
$X\in \L^{1,2}(\calh)$ and that $DX^n\to DX$ weakly
in the space $L^2(\Omega\times [0,T]\times [0,T];\calh)$.
Passing to the limit in (\ref{derivmallappr}) is easily justified
and this proves equation (\ref{derivmall}).
The estimate (\ref{stimasuderivmall})
on $DX$ can be proved in the same way as
it was done for $DX^n$.

We note that for any fixed $s\in [0,T]$, the estimate $|D_\sigma
X^n_s|\le C(s-\sigma)^{\beta-1}$ also shows that $DX_s^n$ is
bounded in the space $L^2(\Omega\times [0,T];\calh)$. Arguing as
before we conclude that $X_s\in \D^{1,2}(\calh)$ for every $s$. The
estimate (\ref{stimasuderivmall}) implies that $DX_s\in L^\infty
(\Omega;L^2([0,T]; \calh))$.

The continuity statement can be proved as in
\cite{DebFuTe}, Proposition 3.4.
\qed

We still set $X_s = X_s^{ 0,x}$, for simplicity.
Given a function $w:[0,T]\times \calh\to\R$,
we investigate the existence of the
joint quadratic variation
of the process $\{w(s,X_s), \,s\in [0,T]\}$
with the Brownian motion $W$ on an interval
$[0,s]\subset [0,T)$. As usual, this is defined
  as the limit in probability of
  $$
\sum_{i=1}^{n}
(w(s_i,X_{s_i}))-w(s_{i-1}, X_{s_{i-1}}))
(W_{s_i}-W_{s_{i-1}})
$$
where $\{s_i\}$,
$0=s_0<s_1<\cdots <s_n=s$ is an
arbitrary subdivision of
$[0,s]$ whose mesh tends to $0$. We do not require that this
convergence takes place uniformly in time.
This definition is easily adapted to an arbitrary interval
of the form $[t,s]\subset [0,T)$.
Existence of the joint quadratic variation is
not trivial. Indeed, due to the occurrence
of convolution type integrals in
the definition of mild solution,
it is not obvious that the process $X$
is a semimartingale. Moreover, even in this case,
 the process $w(\cdot,X)$ might fail to be a
 semimartingale
if $w$ is not twice differentiable, since the It\^o formula does
 not apply. Nevertheless, the following result holds true.
Its proof could be deduced from generalization of
some results obtained in
\cite{NuPa} to the infinite-dimensional case,
but we prefer to give a
simpler direct proof.

\begin{proposition}\label{ucompostox}
Suppose that $w\in C([0,T)\times \calh;
\mathbb{R})$ is
G\^{a}teaux differentiable with respect to $x$,
and that for every $s<T$
there exist constants $K$ and $m$ (possibly
depending on $s$) such that
\begin{equation}\label{stima_u_e_nabla_u}
 |w(t,x)|\leq K(1+|x|)^m,
\qquad |\nabla w (t,x) |\leq K(1+|x|)^m,\qquad t\in [0,s],\;
x\in H.
\end{equation}

Assume that for every
$t\in [0,T)$, $x\in \calh$, $\beta \in(0,\frac{1}{2}+\frac{\theta}{4})$,
 the linear
operator $k\to \nabla w(t,x)(\lambda-A)^{1-\beta}k$ (a priori
defined for $k\in D(\lambda-A)^{1-\beta}$)
has an extension to a bounded linear
operator $\calh\to \mathbb{R}$, that we
 denote by $[\nabla w(\lambda -A)^{1-\beta}](t,x)$.

Moreover assume that
the map $(t,x,k)\to [\nabla w(\lambda-A)^{1-\beta}](t,x)k$
is continuous $[0,T)\times \calh\times \calh \to \mathbb{R}$.

For $t\in [0,T)$, $x\in \calh$, let
$\{X_s^{t,x},\, s\in [t,T]\}$ be the solution of
equation (\ref{eqstatonocontrol}).
Then the process $\{w(s,X_s^{t,x}),\, s\in [t,T]\}$
   admits a
  joint quadratic variation process
 with $W$, on every interval
$[t,s]\subset [t,T)$,
 given by
  $$
  \int_t^s[\nabla w(\lambda -A)^{1-\beta}](r,X_r^{t,x})\,
(\lambda -A)^{\beta}D_\lambda\;dr.
$$
\end{proposition}

{\bf Proof.} For simplicity we take
$t=0$, and we write $X_s=X_s^{0,x}$,
$w_s=w(s,X_s)$. 
It follows from the assumptions that the map
$(t,x,h)\to \nabla w(t,x)h$ is also continuous on
$[0,T)\times \calh\times \calh$.
By the chain rule for
the Malliavin derivative operator (see \cite{FuTe} for details),
 it follows that for every $s <T$
we have $w_s\in\D^{1,2}(\R)$ and
  $Dw_s=\nabla w(s,X_s)DX_s$.

In order to compute the joint quadratic variation of
$w$ and $W$ on a fixed interval $[0,s]\subset [0,T)$.
we take $0=s_0<s_1<\cdots <s_n=s$, a subdivision of
$[0,s]\subset [0,T]$ with mesh $\delta=\max_i (s_i-s_{i-1})$.
By well-known rules of Malliavin
calculus (see \cite{NuPa}, Theorem 3.2, or \cite{GrPa},
Proposition 2.11) we have
$$\begin{array}{l}\dis
(w_{s_i}-w_{s_{i-1}})
(W_{s_i}-W_{s_{i-1}})=
\int^{s_i}_{s_{i-1}}
D_\sigma^j (w_{s_i}-w_{s_{i-1}})\; d\sigma
+\int^{s_i}_{s_{i-1}}
(w_{s_i}-w_{s_{i-1}})\hat{d}W_\sigma^j,
\end{array}
$$
where we use the symbol
$\hat{d}W$ to denote the Skorohod integral.
We note that $D_\sigma w_{s_{i-1}}=0$ for $\sigma>s_{i-1}$.
Therefore setting
$
U_\delta(\sigma)=\sum_{i=1}^{n}
(w_{s_i}-w_{s_{i-1}})\; 1_{(s_{i-1},s_{i}]}(\sigma)$
we obtain
$$
\sum_{i=1}^{n}
(w_{s_i}-w_{s_{i-1}})
(W_{s_i}^j-W_{s_{i-1}}^j) =
\int_0^s U_\delta(\sigma)\; \hat{d}W_\sigma^j
+\sum_{i=1}^{n}\int^{s_i}_{s_{i-1}}
\nabla w(s_i,X_{s_i})\; D_\sigma^jX_{s_i}\;d\sigma.
$$
Recalling (\ref{derivmall}) we obtain
$$\begin{array}{l}\dis
\sum_{i=1}^{n}
(w_{s_i}-w_{s_{i-1}})
(W_{s_i}^j-W_{s_{i-1}}^j) =
\int_0^s U_\delta (\sigma)\; \hat{d}W_\sigma^j+
\sum_{i=1}^{n}\int^{s_i}_{s_{i-1}}
\nabla w(s_i,X_{s_i})\;e^{(s_i-\sigma)A}B\;d\sigma
\\\qquad\qquad\dis
+\sum_{i=1}^{n}\int^{s_i}_{s_{i-1}}
\nabla w(s_i,X_{s_i})\int_\sigma^{s_i}
e^{(s_i-r)A}\nabla F(r,X_r)\;D_\sigma X_r\;dr\;\;d\sigma
=:I_1+I_2+I_3.
\end{array}
$$
Now we let the mesh $\delta$
 tend to $0$.
Following proposition 3.5 \cite{DebFuTe}, we can
prove that $I_1\to 0$ in $L^2(\Omega, \R)$, 
$$I_2\to
  \int_0^s[\nabla w(\lambda -A)^{1-\beta}](r,X_r)\,
(\lambda -A)^{\beta}D_\lambda\;dr,\qquad \P-a.s.
$$
and $I_3\to 0$, $\P$-a.s..
\qed

\section{The backward stochastic differential
equation}

 We consider the following backward stochastic
differential equation:
\begin{equation}\label{BSDE}\left\{\begin{array}{l}
  \dis dY_s^{t,x}=-\Psi(s,X_s^{t,x},Z_s^{t,x})\;ds+Z_s^{t,x} \;dW_s,\qquad s\in [0,T],\\
  Y_T=\Phi(X_T^{t,x}),
\end{array}\right.
\end{equation}
for the unknown real processes  $Y^{t,x}$ and $Z^{t,x}$, also denoted by $Y$ and $Z$.
The equation is understood in the usual way: $\P$-a.s.,
\begin{equation}\label{BSDEmild}
Y_s^{t,x}+\int_s^TZ_r^{t,x} \;dW_r=
\Phi(X_T^{t,x})+\int_s^T\Psi(r,X_r^{t,x},Z_r^{t,x})\;dr,\qquad s\in [0,T],
\end{equation}
but we will use the shortened notation above for equation (\ref{BSDE})
and similar equations to follow.
In (\ref{BSDE}) and (\ref{BSDEmild}),
$t\in [0,T]$ and $x\in \calh$ are given and the process
$X^{t,x}$ is the solution of (\ref{eqstatonocontrol}), with the
convention that $X_s^{t,x}=x$ for $s\in [0,t)$.
On the generator $\Psi$ and on the final datum $\Phi$ we make the following assumptions:
\begin{hypothesis}\label{ipotesisufipsi_backward}

\item[1)] $|\Phi(x_1)-\Phi(x_2)|_\calh \leq
C_{\Phi}(1+|x_1|+|x_2|)|x_2-x_1|$ for all $x_1$, $x_2$ in $\calh$.
 \item[2)] There
exists a constant $C_{\psi}$ such that
 $|\Psi(t,x_1,z)-\Psi(t,x_2,z)|\leq C_{\psi}(1+|x_1|+|x_2|)|x_2-x_1|$
 for all $x_1$, $x_2$ in $\calh$, $z\in \mathbb{R}$ and $t\in [0,T]$ and
  $|\Psi(s,x,z_1)-\Psi(s,x,z_2)|\le C_\psi |z_1-z_2|, $ for every
$s\in [0,T]$, $x\in \calh$, $z_1,z_2\in\R$.
\item[3)]
$\sup_{s\in[0,T]} |\Psi(s,0,0)| \leq C_{\ell}.$
 \item[4)] $\Phi\in \calg^1(\calh)$ and for almost every $s\in [0,T]$ the map $\Psi(s,\cdot,\cdot)$ is
G\^{a}teaux differentiable on $\calh\times \R$ and the maps $(x,h,z)\to
\nabla_x\Psi(s,x,z)h$ and $(x,z,\zeta)\to
\nabla_z\Psi(s,x,z)\zeta$ are continuous on
$\calh\times \calh\times \R$ and $\calh\times \R\times \R$ respectively.
\end{hypothesis}

\begin{proposition}\label{bsderegolare}
 \begin{description}\label{esunYZ}
    \item[$1)$] For all $x\in \calh$, $t\in [0,T]$ and $p\in [2,\infty)$
    there exists
a unique pair of processes $(Y^{t,x},Z^{t,x})$ with $Y^{t,x}\in
L^p_{\mathcal{P}}(\Omega,C([0,T],\mathbb{R}))$, $Z^{t,x}\in
L^p_{\mathcal{P}}(\Omega,L^2([0,T],\mathbb{R}))$ solving
(\ref{BSDE}); in the following we denote such a solution by
$(Y_\cdot^{t,x},Z_\cdot^{t,x})$.
    \item[$2)$] The map $(t,x)\to (Y_\cdot^{t,x},Z_\cdot^{t,x})$ is
    continuous from $[0,T]\times \calh$ to
    $L^p_{\mathcal{P}}(\Omega,C([0,T],\mathbb{R}))\times
    L^p_{\mathcal{P}}(\Omega,L^2([0,T],\mathbb{R}))$.
    \item[$3)$] For all $t\in [0,T]$ the map
    $x\to (Y_\cdot^{t,x},Z_\cdot^{t,x})$
    is G\^ateaux differentiable as a map from $\calh$ to
    $L^p_{\mathcal{P}}(\Omega,C([0,T],\mathbb{R}))\times
    L^p_{\mathcal{P}}(\Omega,L^2([0,T],\mathbb{R}))$;
    moreover the map
    $(t,x,h)\to (\nabla_x Y_\cdot^{t,x}h$, $\nabla _x Z_\cdot^{t,x}h)$ is
    continuous from $[0,T]\!\times\! \calh\!\times\! \calh$ to
    $L^p_{\mathcal{P}}(\Omega,C([0,T],\mathbb{R}))\times
L^p_{\mathcal{P}}(\Omega,L^2([0,T],\mathbb{R}))$.
 \item[$4)$]
The following equation holds for all $t\in [0,T]$, $x,h \in \calh$.
\begin{equation}\label{BSDE-deriv}\left\{\begin{array}{l}
  \dis d\;\nabla_x Y_s^{t,x}h=-\nabla_x\Psi(s,X_s^{t,x},Z_s^{t,x})
  \,\nabla_x X_s^{t,x}h\; ds\\ \\
 \qquad\quad\qquad \dis -\nabla_z \Psi(s,X_s^{t,x},Z_s^{t,x})
 \,\nabla_x Z_s^{t,x}h\;ds+ \nabla_x Z_s^{t,x}h\;
 dW_s,\\ \\
  \dis\nabla_xY_T^{t,x}h=\nabla_x\Phi(X_T^{t,x})\,\nabla_xX_T^{t,x}h,
  \qquad \qquad s\in [t,T].
\end{array}\right.
\end{equation}
\end{description}
\end{proposition}
{\bf Proof.} The claim follows directly from Proposition 4.8 in
\cite{FuTe}, from Proposition \ref{regequazforward} above and from
the chain rule (in the form stated in Lemma 2.1 of \cite{FuTe}). \qed

\begin{remark}{\em
The inequality (\ref{stimaes}), for $\alpha=0$, together with the inequality
(4.9) in \cite{FuTe}, implies that there exists a constant
$C_{Y,p}$ such that for all $t\in [0,T]$ and $x\in \calh$
$$\E\sup_{s\in [0,T]}|Y_s^{t,x}|^p+\E\left(\int_0^T |
Z_s^{t,x}|^2 ds\right)^{p/2} \leq C_{Y,p} (1+|x|)^{2p}. $$}\end{remark}
\begin{remark}{\em $Y_t^{t,x}$ is adapted both to the
$\sigma$-field $\sigma\{W_s: s\in [0,t]\}$ and to the
$\sigma$-field $\sigma\{W_s-W_t: s\in [t,T]\}$. Thus $Y_t^{t,x}$
and $\nabla Y_t^{t,x}h$, $x,h\in H$ are deterministic.
}\end{remark}

 Proposition \ref{regequazforward2} yields the
following further regularity result. 

\begin{proposition}\label{regback}
For every $\alpha \in [0,1/2)$, $p\in [2,\infty)$ there exist two
families of processes $$\left\{P^{\alpha}(s,t,x)k : s\in
[0,T]\right\} \hbox{ and }\left\{ Q^{\alpha}(s,t,x)k:s\in
[0,T]\right\}; \quad \hbox{$t\in [0,T)$, $x\in \calh$, $k\in \calh$}$$
with $P^{\alpha}(\cdot,t,x)k\in
L^p_{\mathcal{P}}(\Omega,C([0,T],\mathbb{R}))\hbox{ and }
Q^{\alpha}(\cdot,t,x)k )\in
L^p_{\mathcal{P}}(\Omega,L^2([0,T],\mathbb{R}))$ such that if
$k\in D(\lambda-A)^{\alpha}$, $t\in [0,T)$, $x\in \calh$, then $\mathbb{P}$-a.s.
\begin{equation}\label{def-di-P}
    P^{\alpha}(s,t,x)k=\left\{\begin{array}{ll}
  \nabla_x Y_s^{t,x}(\lambda-A)^{\alpha}k & \text{ for all }s \in [t,T],  \\
  \nabla_x Y_t^{t,x}(\lambda-A)^{\alpha}k& \text{ for all }s \in [0,t), 
\end{array} \right.
\end{equation}
\begin{equation}\label{def-di-Q}Q^{\alpha}(s,t,x)k=\left\{\begin{array}{ll}
  \nabla_x Z_s^{t,x}(\lambda-A)^{\alpha}k & \hbox{ for a.e. $s \in [t,T]$,} \\
  0& \hbox{ if $s \in [0,t)$}.
\end{array} \right.\end{equation}
Moreover the map $ (t,x,k)\to P^{\alpha}(\cdot,t,x)k$ 
and the map $(t,x,k)\to Q^{\alpha}(\cdot,t,x)k $
are continuous from $[0,T)\times \calh\times \calh$ to
$L^p_{\mathcal{P}}(\Omega,C([0,T],\mathbb{R}))$ and 
linear with respect to $k$.

Finally there exists a constant $C_{\nabla Y,\alpha,p}$ such that
\begin{equation}\label{stimadiP}
    \E\sup_{s\in [0,T]}|P^{\alpha}(s,t,x)k|_\calh^p+\E\left(\int_0^T
|Q^{\alpha} (s,t,x)k|_{(\mathbb{R})} ds\right)^{p/2} \leq
C_{\nabla Y,\alpha,p} (T-t)^{-\alpha p} (1+|x|_\calh)^p|k|_\calh^p.
\end{equation}
\end{proposition}
{\bf Proof.} Let, for $t\in [0,T]$, $x\in \calh$, $k\in D(\lambda-A)^{\alpha}$,
$P^{\alpha}(\cdot,t,x)k$ and $Q^{\alpha}(\cdot,t,x)k $ be defined by
(\ref{def-di-P}) and (\ref{def-di-Q}) respectively.

\noindent By Proposition \ref{bsderegolare}
  the map $k\to(P^{\alpha}(\cdot,t,x)k,Q^{\alpha}(\cdot,t,x)k )$ is
 a bounded linear operator from $D(\lambda-A)^{\alpha}$ to
$L^p_{\mathcal{P}}(\Omega,C([0,T],\mathbb{R}))\times
L^p_{\mathcal{P}}(\Omega,L^2([0,T],\mathbb{R}))$. Moreover $
(P^{\alpha}(\cdot,t,x)k,Q^{\alpha}(\cdot,t,x)k) $ solves the equation
\begin{equation}\label{equaz_P}\left\{\begin{array}{l}
  \dis d P^{\alpha}(s,t,x)k=-1_{[t,T]}(s) \nabla_x\Psi(s,X_s^{t,x},Z^{t,x})
  \,\nabla_x X_s^{t,x}(\lambda-A)^{\alpha }k \;ds\\ \\
 \qquad\quad\qquad \dis
 -\nabla_z \Psi(s,X_s^{t,x},Z_s^{t,x})\;Q^{\alpha}(s,t,x)k\; ds
 + Q^{\alpha}(s,t,x)k\; dW_s,\\ \\
  \dis P^{\alpha}(T,t,x)k=\nabla_x\Phi(X_s^{t,x})\,
  \nabla_xX_T^{t,x}(\lambda-A)^{\alpha }k, \qquad \qquad s\in [t,T].
\end{array}\right.
\end{equation}
By (\ref{defdiTheta}) equation (\ref{equaz_P}) can be rewritten
\begin{equation}\label{equaz_P_2}
    \left\{\begin{array}{l}
  \dis \!d P^{\alpha}(s,t,x)k=\nu(s,t,x)k\;
   ds -\nabla_z \Psi(s,X_s^{t,x},Z_s^{t,x})Q^{\alpha}(s,t,x)k ds+
   Q^{\alpha}(s,t,x)k\; dW_s\\ \\
  \dis \!P^{\alpha}(T,t,x)k= \eta(t,x)k,
   \qquad \qquad s\in [0,T],
\end{array}\right.
\end{equation}
where
$$\begin{array}{l}
  \dis \nu(s,t,x)k= -1_{[t,T]}(s)
  \nabla_x\Psi(s,X_s^{t,x},Z_s^{t,x})
  \left(\Theta^{\alpha} (s,t,x)k+ e^{(s-t)A} (\lambda-A)^{\alpha}k\right),
  \\
  \eta(t,x)k=
  \nabla_x\Phi(X_T^{t,x})\left(\Theta^{\alpha}(T,t,x)k
  + e^{(T-t)A}(\lambda-A)^{\alpha}k\right).
\end{array}
$$
Now we choose arbitrary $k\in \calh$ and notice that $\nu(s,t,x)k$ and
$\eta(t,x)k$ can still be defined by the above formulae. Remark
\ref{stima_nabla_Psi}, and relations (\ref{stimaes}), with $\alpha=0$,
(\ref{stimadi nabla theta}) yield:
$$
\begin{array}{rcl}
  \dis \E\left(\int_0^T |\nu (s,t,x)k|^2ds \right)^{p/2}&\!\!\!\!\!\!\!\!\!\leq\!\!\!& \dis c_1
\E\!\left(\int_t^T (1+ |X_s^{t,x}|)^2\left( |\Theta^{\alpha}
(s,t,x)k| +(s-t)^{-\alpha}|k|\right)^2ds\right)^{p/2} \\
  &\!\!\!\!\!\!\!\!\!\!\!\! \leq \!\!\! &
  \dis c_2 \left[(T-t)^{p/2}+(T-t)^{(1-2\alpha)p/2}\right] (1+|x|)^p|k|^p
  \leq c_3(1+|x|)^p|k|^p,
\end{array}
$$
where $c_1$, $c_2$ and $c_3$ are suitable constants independent on
$t,x,k$. In the same way
$$
\begin{array}{rcl}
  \dis \E |\eta (t,x)k|^p&\leq& \dis c_4
\E\Big( (1+ |X_T^{t,x}|)\left( |\Theta^{\alpha}
(T,t,x)k| +(T-t)^{-\alpha}|k|\right)\Big)^{p} \\
  & \leq & \dis c_5 (T-t)^{-p\alpha}(1+|x|)^p|k|^p.
\end{array}
$$
By Proposition 4.3 in \cite{FuTe}, for all $k\in \calh$ there
exists a unique pair
$(P^{\alpha}(\cdot,t,x)k,Q^{\alpha}(\cdot,t,x)k )$
 belonging to
$L^p_{\mathcal{P}}(\Omega,C([0,T],\mathbb{R}))\times
L^p_{\mathcal{P}}(\Omega,L^2([0,T],\mathbb{R}))$ and solving
equation (\ref{equaz_P_2}); moreover (\ref{stimadiP}) holds. The map
$k\to (P^{\alpha}(\cdot,t,x)k,Q^{\alpha}(\cdot,t,x)k )$ is clearly
linear, so we can conclude that the required extension exists. The
proof of its continuity can be achieved as in \cite{DebFuTe}, proposition 4.4.
\qed

\begin{corollary}\label{proregu} Setting $v(t,x)=Y_t^{t,x}$, we have
$v\in C([0,T]\times \calh;
\mathbb{R})$ and there exists a constant $C$ such that
$|v(t,x)|\leq C\, (1+|x|)^2$, $t\in [0,T]$, $x\in \calh$.
 Moreover $v$ is
G\^{a}teaux differentiable with respect to $x$
on $[0,T]\times \calh$ and the map
$(t,x,h)\to \nabla v(t,x)h$ is continuous.

For all $\alpha\in [0,1/2)$, $t\in [0,T)$ and $x\in \calh$ the linear
operator $k\to \nabla v(t,x)(\lambda-A)^{\alpha}k$ - a priori
defined for $k\in D(\lambda-A)^{\alpha}$ - has an extension to a bounded linear
operator $\calh\to \mathbb{R}$, that we denote by
$[\nabla v(\lambda-A)^{\alpha}](t,x)$.

Finally the map $(t,x,k)\to [\nabla v(\lambda-A)^{\alpha}](t,x)k$
is continuous $[0,T)\times \calh\times \calh \to \mathbb{R}$ and there
exists $C_{\nabla v,\alpha}$ for which:
\begin{equation}\label{stima_nabla_u}
|[\nabla v (\lambda-A)^{\alpha}](t,x)k |\leq C_{\nabla v,\alpha}
(T-t)^{-\alpha}(1+|x|_\calh)|k|_\calh, \qquad t\in [0,T),\,\,
x, k\in \calh.
\end{equation}
\end{corollary}
{\bf Proof.}
We recall that
  $Y_t^{t,x}$ is
  deterministic.
Since the map $(t,x)\rightarrow Y^{t,x}$ is
  continuous with values in
  $L^p_{\mathcal{P}}(\Omega,C([0,T],\mathbb{R}))$, $p\geq 2$,
  then the map $(t,x)\rightarrow Y_t^{t,x}$ is
  continuous with values in
  $L^p(\Omega,\mathbb{R})$ and so the map $(t,x)\rightarrow 
  Y_t^{t,x}=v(t,x)$ is continuous with values in $\mathbb{R}$.

Similarly, $ \nabla_x v(t,x)= \nabla_x Y_t^{t,x}$ exists and
has the required continuity properties, by Proposition
\ref{esunYZ}.

Next we notice that
$P^{\alpha}(t,t,x)k=
  \nabla_x Y_t^{t,x}(\lambda-A)^{\alpha}k$. The existence
  of the required extensions
and its continuity are direct consequences of Proposition
\ref{regback}. Finally
the estimate (\ref{stima_nabla_u}) follows from
(\ref{stimadiP}).
\qed
\begin{remark}{\em It is evident by construction
that the law of $Y^{t,x}$ and consequently the function $v$
depends on the law of the Wiener process $ W$ but not on the
particular probability $\mathbb{P}$ and Wiener process $ W$ we
have chosen. }\end{remark}

\begin{corollary}\label{identmarkov} For every $t\in [0,T]$, $x\in H$
we have, $\P$-a.s.,
\begin{equation}\label{identifY}
Y_s^{t,x}=v(s ,X_s^{t,x}),\qquad
{\rm \; for\;all\; }s \in [ t,T],
\end{equation}
\begin{equation}\label{identifZ}
    Z_s^{t,x}=
[\nabla v(\lambda-A)^{1-\beta}](s ,X_s^{t,x})\;(\lambda-A)^{\beta}D_\lambda,\qquad
{\rm \; for\;almost\;all\;}s \in [ t,T].
\end{equation}
\end{corollary}

{\bf Proof.}
We start from the well-known equality: for $0\le t\le r\le T$,
$\P$-a.s.,
$$
X_s^{t,x}=X_s^{r,X_r^{t,x}},\qquad
{\rm \; for\;all\; }s \in [ r,T].
$$
It follows easily from the uniqueness of the backward equation
(\ref{BSDE}) that
$\P$-a.s.,
$$
Y_s^{t,x}=Y_s ^{r,X_r^{t,x}},\qquad
{\rm \; for\;all\; }s \in [ r,T].
$$
Setting $s=r$ we arrive at (\ref{identifY}).

To prove (\ref{identifZ}) we note that it follows
immediately from the backward equation
(\ref{BSDE}) that the joint quadratic variation of
$\{Y_s^{t,x},\;s \in [ t,T]\}$ and $W$
on an arbitrary interval $[t,s]\subset [t,T)$ is equal
to $\int_t^s{Z}^j \; dr$. By
(\ref{identifY}) the same result can be obtained
by considering the joint quadratic variation of
$\{v(s,X_s^{t,x}),\;s \in [ t,T]\}$ and $W$. An application
of Proposition \ref{ucompostox} (whose assumptions hold
true by Corollary \ref{proregu}) leads to the identity
$$
\int_t^s{Z}_r \; dr=
  \int_t^s[\nabla v(\lambda -A)^{1-\beta}](r,X_r^{t,x})\,
(\lambda -A)^{\beta}D_\lambda\;dr,
$$
and (\ref{identifZ}) is proved.
\qed

\section{The Hamilton-Jacobi-Bellman equation}

In this section the aim is to solve a second order partial differential equation,
where the second order differential operator is the generator of the Markov process
$\{X_s^{t,x},s\in[t,T]\}$, solution of equation (\ref{eqstatonocontrol}).
Namely we are interested in Hamilton Jacobi Bellman equations
related to a control problem that we present in the next section.

Let us consider again the solution $X_s^{t,x}$ of equation
(\ref{eqstatonocontrol}) and denote by $P_{t,s}$
its transition semigroup:
$$
P_{t,s}[\phi](x)=\E\, \phi(X_s^{t,x}),\qquad x\in \calh,\;0\le t\le s\le T,
$$
for any bounded measurable $\phi:\calh\to \R$.
We note that by the estimate (\ref{stimaes}) (with $\alpha=0$)
this formula
is meaningful for every $\phi$ with polynomial growth.
In the following $P_{t,s}$ will be considered
as an operator acting on this class of functions.

Let us denote by $\call_t$ the generator of $P_{t,s}$,
formally:
$$
\call_t[\phi](x)=\frac{1}{2}
\< \nabla^2\phi(x)B,B\>
+ \< Ax+F(t,x),\nabla\phi(x)\>,
$$
where $\nabla\phi(x)$ and $\nabla^2\phi(x)$
are first and second G\^ateaux derivatives of
$\phi$ at the point $x\in \calh$ (here they are
identified with elements of $\calh$ and $L(\calh)$ respectively).

The Hamilton-Jacobi-Bellman equation for the optimal
control problem is
\begin{equation}\label{HJBformale}
  \left\{\begin{array}{l}\dis
\frac{\partial v(t,x)}{\partial t}+\call_t [v(t,\cdot)](x) =
-\Psi (t, x,\nabla v(t,x)B),\qquad t\in [0,T],\,
x\in \calh,\\
\dis v(T,x)=\Phi(x).
\end{array}\right.
\end{equation}
This is a nonlinear parabolic equation
for the unknown function $v:[0,T]\times \calh\to\R$.
The operators $\call_t$ are very degenerate, since
the space $\calh$ is infinite-dimensional but
the noise $W$ is a real Wiener process.

Now we consider the variation of
constants formula for
 (\ref{HJBformale}):
$$
  v(t,x) =P_{t,T}[\Phi](x)-\int_t^TP_{t,s}[
\Psi (s, \cdot,
\nabla v(s,\cdot)B
](x)\; ds,\qquad t\in [0,T],\,
x\in \calh,
$$
where we remember $B=(\lambda-A)D_\lambda$.
This equality is still formal, since the
term $(\lambda-A)D_\lambda$ is not defined. However with a slightly
different interpretation we arrive at the following precise
definition:
\begin{definition}\label{defdisoluzionemild}
Let $\beta\in[0,\frac{1}{2})$. We say that a function
$v:[0,T]\times \calh\to\R$ is a mild solution of the
Hamilton-Jacobi-Bellman equation
(\ref{HJBformale}) if the following conditions hold:
\begin{enumerate}
  \item[(i)]
$v\in C([0,T]\times \calh;
\mathbb{R})$ and there exist constants $C,m\ge 0$ such that
$|v(t,x)|\leq C\, (1+|x|)^m$, $t\in [0,T]$, $x\in \calh$.

 \item[(ii)] $v$ is
G\^{a}teaux differentiable with respect to $x$
on $[0,T)\times \calh$ and the map
$(t,x,h)\to \nabla v(t,x)h$ is continuous
$[0,T)\times \calh\times \calh\to \R$.

 \item[(iii)]
For all $t\in [0,T)$ and $x\in \calh$ the linear
operator $k\to \nabla v(t,x)(\lambda-A)^{1-\beta}k$ (a priori
defined for $k\in D(\lambda-A)^{1-\beta}$) has an extension to a bounded linear
operator $\calh\to \mathbb{R}$, that we denote
by $[\nabla v(\lambda-A)^{1-\beta}](t,x)$.

Moreover the map $(t,x,k)\to [\nabla v(\lambda-A)^{1-\beta}](t,x)k$
is continuous $[0,T)\times \calh\times \calh \to \mathbb{R}$ and there
exist constants $C,m\ge0$, $\kappa\in [0,1)$ such that
\begin{equation}\label{stima_grad_u_in_generale}
|[\nabla v (\lambda-A)^{1-\beta}](t,x) |_{\calh^*}\leq C
(T-t)^{-\kappa}(1+|x|)^m, \qquad t\in [0,T),\,\,
x\in \calh.
\end{equation}

  \item[(iv)] the following equality holds for
  every $t\in [0,T]$, $x\in \calh$:
  \begin{equation}\label{solmild}
  v(t,x) =P_{t,T}[\Phi](x)+\int_t^TP_{t,s}\left[
\Psi \Big(s, \cdot,
[\nabla v(\lambda-A)^{1-\beta}](s,\cdot)\;(\lambda-A)^{\beta}D_\lambda\Big)
\right](x)\; ds.
\end{equation}

\end{enumerate}
\end{definition}

We assume that $\Phi$ and $\Psi$ satisfy hypotheses 
\ref{ipotesisufipsi_backward} and using the estimate
(\ref{stimaes}) (with $\alpha=0$)
it is easy to conclude that formula (\ref{solmild}) is meaningful.

\begin{theorem}\label{main}
Assume Hypotheses \ref{ipotesiconcrete},
\ref{ipotesiconcretecosto} and \ref{ipotesisufipsi_backward} then there
exists a unique mild solution of the Hamilton-Jacobi-Bellman
equation (\ref{HJBformale}).
The solution $v$ is given by the formula
$$
    v(t,x) =Y_t^{t,x},
  $$
where $(X,Y,Z)$ is the solution of the forward-backward
system (\ref{eqstatonocontrol})-(\ref{BSDEmild}).

\end{theorem}

{\bf Proof.}
{\em Existence.}
By Corollary \ref{proregu} the solution $v$ has the
regularity properties stated in Definition
\ref{defdisoluzionemild}.
In order to verify that
equality (\ref{solmild}) holds
we first fix $t\in [0,T]$ and
$x\in \calh$ and note that
the backward equation (\ref{BSDE}) gives
$$
 Y_t^{t,x}+\int_t^TZ_s^{t,x} \; dW_s =
  \Phi(X_T^{t,x})
  +\int_t^T\Psi \Big(s ,X_s^{t,x},Z_s^{t,x}\Big)\; ds .
$$
Taking expectation we obtain
\begin{equation}\label{quasimild}
 v(t,x)=
  P_{t,T}[\Phi](x)
  +\E\int_t^T\Psi \Big(s ,X_s^{t,x}, Z_s^{t,x}
  \Big)\; ds .
\end{equation}
Now we recall that by Corollary \ref{identmarkov} we have
$$
    Z_s^{t,x}=
[\nabla v(\lambda-A)^{1-\beta}](s ,X_s^{t,x})\;(\lambda-A)^{\beta}D_\lambda,\qquad
\P {\rm -a.s. \; for\;a.a.\;}s \in [ t,T].
  $$
It follows that
$$
\E\int_t^T
\Psi \Big(s , X_s^{t,x}
Z_s^{t,x}\Big)
\; ds=\int_t^TP_{t,s }\left[
\Psi \Big(s , \cdot,
[\nabla v(\lambda-A)^{1-\beta}](s ,\cdot)\;(\lambda-A)^{\beta}D_\lambda
\Big)
\right](x)\; ds.
$$
Comparing with (\ref{quasimild}) gives the required
equality (\ref{solmild}).

{\em Uniqueness.}
Let $v$ be a mild solution. We fix $t\in [0,T]$
and $x\in \calh$ and look for a convenient
expression for the process $v(s,X_s^{t,x})$, $s\in [t,T]$. 
By ``standard'' arguments (see e.g. \cite{FuTe}), by the Markov property 
of $X$ the process
$v(s,X_s^{t,x})$, $s\in [t,T]$ is a (real) continuous
semimartingale, and, by the representation theorem for martingales, there exists
$\widetilde{Z}\in L^2_\calp(\Omega \times [t,T]; \R)$ such that
its canonical decomposition into its continuous martingale 
part and its continuous finite variation part is given by
\begin{equation}\label{scomposizione}
\begin{array}{l}\dis
  v(s,X_s^{t,x})
=v(t,x)+
\int_t^s\widetilde{Z}_r\; dW_r
\\\dis\qquad\qquad
+
\int_t^s
\Psi \Big(r, X_r^{t,x},
[\nabla v(\lambda-A)^{1-\beta}](r,X_r^{t,x})\;(\lambda-A)^{\beta}D_\lambda
\Big)\; dr.
\end{array}
\end{equation}

By computing the joint quadratic variations 
of both sides of (\ref{scomposizione})
we have $\P$-a.s.,
$
[\nabla v(\lambda -A)^{1-\beta}](s,X_s^{t,x})\,
(\lambda -A)^{\beta}D_\lambda=
\widetilde{Z}_s,
$ so
substituting into (\ref{scomposizione}) and taking into account
that $v(T,X_T^{t,x}))=\Phi(X_T^{t,x}))$ we obtain,
for $s\in [t,T]$,
$$
\begin{array}{l}\dis
  v(s,X_s^{t,x})
+
\int_s^T
[\nabla v(\lambda -A)^{1-\beta}](r,X_r^{t,x})\,
(\lambda -A)^{\beta}D_\lambda\;dW_r
\\\dis\qquad\qquad
 =\Phi(X_T^{t,x})
+\int_s^T
\Psi \Big(r, X_r^{t,x},
[\nabla v(\lambda-A)^{1-\beta}](r,X_r^{t,x})\;(\lambda-A)^{\beta}D_\lambda
\Big)\; dr.
\end{array}
$$
 Comparing with
the backward equation (\ref{BSDE})
we note that the pairs
$$
\Big(Y_s^{t,x}, Z_s^{t,x}\Big)\;{\rm and}\;
\Big(v(s,X_s^{t,x}), [\nabla v(\lambda -A)^{1-\beta}](s,X_s^{t,x})\,
(\lambda -A)^{\beta}D_\lambda\Big), \; s\in [t,T],
$$
solve the same equation. By uniqueness, we have
$Y_s^{t,x}=v(s,X_s^{t,x})$, $s\in [t,T]$, and
setting $s=t$ we obtain $Y_t^{t,x}=v(t,x)$.
\qed

\section{Synthesis of the optimal control}

At first we introduce a ``concrete" cost functional: 
let $y(s,\xi )$ solution of equation (\ref{eqconcreta1}).
Let us consider the following cost functional
\begin{equation}\label{costoconcreto1}
J(t,x,u(\cdot))=\E \int_t^T\int_0^{+\infty} \ell(s,\xi,
y(s,\xi),u(s))\;d\xi\;ds +\E \int_0^{+\infty} \phi(\xi,
y(T,\xi))\;d\xi.
\end{equation}
In this section we assume that the following holds:
\begin{hypothesis}\label{ipotesiconcretecosto}
 $\ell:[0,T]\times[0,+\infty)\times\R\times\calu
\to\R$ and
$\phi:[0,+\infty)\times\R \to\R$ are measurable. Let $\rho(\xi)=\xi^{1+\theta}$,
or $\rho(\xi)=\xi^{1+\theta}\wedge 1$, depending on what weight we are considering to
define the space $\calh$. Assume also:
\begin{enumerate}
 \item[1)] there exist two constant $C_1, C_2$ such that,
for some $\epsilon >0$, for every $\xi\in [0,+\infty)$,
   $y_1,y_2\in \mathbb{R}$
 $$|\phi(\xi, y_1)- \phi(\xi, y_2)| \leq C_1
 \dfrac{\sqrt{\rho(\xi)}}{(1+\xi)^{1/2+\epsilon}}\; |y_1-y_2|
\,+\, C_2 \;\rho(\xi)(|y_1|+|y_2|)\;|y_1-y_2|,
  $$
 moreover $\dis\int_0^{+\infty} |\phi(\xi,0)|d\xi<\infty$;
  \item[2)]
  for every $t\in [0,T]$ and $\xi \in [0,+\infty)$,
  $\ell(t,\xi,\cdot,\cdot):\mathbb{R}^2 \rightarrow \mathbb{R}$
  is continuous.
   Moreover there exists  two constant $C_1, C_2$ such that,
for some $\epsilon >0$, for every $t\in [0,T]$, $\xi\in [0,+\infty)$,
   $y_1,y_2\in \mathbb{R}$, $u\in \calu$,
$$
|\ell(t,\xi, y_1,u)- \ell(t,\xi, y_2,u)| \leq C_1
 \dfrac{\sqrt{\rho(\xi)}}{(1+\xi)^{1/2+\epsilon}}\; |y_1-y_2|
\,+\, C_2\; \rho(\xi)(|y_1|+|y_2|)\;|y_1-y_2|,
 $$
 and for every $t\in [0,T]$
 $$\int_0^{+\infty} \sup_{u\in \calu}
 |\ell(t,\xi,0,u)|\; d\xi \leq C_{\ell}.
  $$
\end{enumerate}
\end{hypothesis}

We notice that in Hypothesis \ref{ipotesiconcretecosto}, the presence 
of the weight $\rho(\xi)$ is natural since we are considering as 
state space the weighted space $\calh$, as well as the presence of
the square integrable function $\dfrac{1}{(1+\xi^{1/2+\epsilon})}$
since $[0,+\infty)$ is not of finite measure with any weight $\rho(\xi)$.

\noindent Further assumptions will be made on the cost functional after 
the following reformulation: we define
$$
L(s,x,u)=\int_0^{+\infty} \ell(s,\xi, x(\xi),u)\;d\xi, \qquad
 \Phi(x)= \int_0^{+\infty} \phi(\xi, x(\xi))\;d\xi,
$$
for $s\in [0,T]$, $x=x(\cdot)\in \calh$,
$u\in \calu$.
The functions $L:[0,T]\times \calh
\times \calu\to\R$ and $\Phi:\calh \to\R$
 are well defined and measurable.
 The cost functional (\ref{costoconcreto1})
can be written in the form
\begin{equation}\label{costoastratto}
J(t,x,u(\cdot))=\E \int_t^T L(s, X_s^u,u_s)\;ds +\E \,\Phi(X^u_T).
\end{equation}
It is easy to show that the cost is finite for any admissible control
$u(\cdot)$.
Moreover for $s\in [0,T]$, $x\in \calh$, $z\in\R$ we define the
hamiltonian as
$$
\Psi(s,x,z)=\inf_{u\in\calu} \{zu+L(s,x ,u) \}.
$$
Since, as it is easy to check, for all $s\in [0,T]$ and all $x\in
\calh$, $ L(s,x,\cdot)$ is continuous on the compact set $\calu$ the
above infimum is attained. Therefore if we define
\begin{equation}\label{defdigammagrande}
\Gamma(s,x,z)=\left\{ u\in \calu: zu+L(s,x ,u)= \Psi(s,x,z)\right\}
\end{equation}
 then
$\Gamma(s,x,z) \neq \emptyset$ for every $s\in [0,T]$, every $x\in
\calh$ and every $z\in \R$. By \cite{AubFr}, see Theorems 8.2.10 and
8.2.11, $\Gamma$ admits a measurable selection, i.e. there exists
a measurable function $\gamma: [0,T]\times \calh \times \R 
\rightarrow \calu$ with $\gamma(s,x,z)\in \Gamma(s,x,z)$ for every
$s\in [0,T]$, every $x\in \calh$ and every $z\in \R$.

\begin{proposition}\label{ipotesiastrattecosto}
Under Hypothesis \ref{ipotesiconcretecosto} the following holds.
\begin{enumerate}
\item[1)] $|\Phi(x_1)-\Phi(x_2)|_\calh \leq
C_{\phi}(1+|x_1|+|x_2|)|x_2-x_1|$ for all $x_1$, $x_2$ in $\calh$.
 \item[2)] There
exists a constant $C_{\psi}$ such that
 $|\Psi(t,x_1,z)-\Psi(t,x_2,z)|\leq C_{\psi}(1+|x_1|+|x_2|)|x_2-x_1|$
 for all $x_1$, $x_2$ in $\calh$, $z\in \mathbb{R}$ and $t\in [0,T]$.
  \item[3)] Setting $C_\calu=\sup\{|u|\,:\,u\in\calu\}$ we have
  $|\Psi(s,x,z_1)-\Psi(s,x,z_2)|\le C_\calu\, |z_1-z_2|, $ for every
$s\in [0,T]$, $x\in \calh$, $z_1,z_2\in\R$.
\item[4)]
$\sup_{s\in[0,T]} |\Psi(s,0,0)| \leq C_{\ell}.$
\end{enumerate}
\end{proposition}


Some of our results are based on the following
assumptions:

\begin{hypothesis}\label{ipotesisupsi}
For almost every $\xi\in [0,+\infty)$ the map $\phi(\xi,\cdot)$
is continuously differentiable on $\R$.
For almost every $s\in [0,T]$ the map $\Psi(s,\cdot,\cdot)$ is
G\^{a}teaux differentiable on $\calh\times \R$
and the maps $(x,h,z)\to
\nabla_x\Psi(s,x,z)h$ and $(x,z,\zeta)\to
\nabla_z\Psi(s,x,z)\zeta$ are continuous on
$\calh\times \calh\times \R$ and $\calh\times \R\times \R$
respectively.
\end{hypothesis}

\noindent From this assumption and from
Hypothesis \ref{ipotesiconcretecosto}
it follows
easily that $\Phi$ is
G\^{a}teaux differentiable on $\calh$ and the map $(x,h)\to
\nabla\Phi(x)h$ is continuous on $\calh\times \calh$. Moreover it
follows that $\Phi$ and $\Psi$ satisfy hypothesis \ref{ipotesisufipsi_backward}.

\begin{remark}\label{stima_nabla_Psi}
{\em From Proposition \ref{ipotesiastrattecosto}
we immediately deduce the following estimates:
$$|\nabla\Phi(x)h|\leq C_{\phi}(1+2|x|)|h|,\quad
|\nabla_x\Psi(t,x,z)h|\leq C_{\psi}(1+2|x|)|h|,\quad
|\nabla_z\Psi(s,x,z)\zeta|\le C_\calu\, |\zeta|.
$$
 }\end{remark}

Hypothesis \ref{ipotesisupsi} involves conditions
on the function $\Psi$, and not on the function  
$\ell$ that determines $\Psi$.
However, Hypothesis \ref{ipotesisupsi} can be
verified in concrete situations, see e.g. example 2.7.1 in
\cite{DebFuTe}.

%
%

The optimal control problem in its strong formulation is to
minimize, for arbitrary $t\in [0,T]$ and $x\in \calh$,
the cost (\ref{costoastratto}),
over all
admissible controls, where
$\{X^u_s\,:\, s\in [t,T]\}$ solves $\P$-a.s.
\begin{equation}\label{eqstatodue}
  \begin{array}{lll}
  \dis
X^u_s&=&\dis e^{(s-t)A}x+\int_t^s e^{(s-r)A}F(r,X^u_r)\;dr +
\int_t^s(\lambda -A)^{1-\beta}e^{(s-r)A}(\lambda
-A)^{\beta}D_\lambda\;dW_r
\\
&+&\dis
\int_t^s(\lambda -A)^{1-\beta}e^{(s-r)A}(\lambda -A)^{\beta}D_\lambda\;u_r\;dr,
\qquad s\in [t,T].
\end{array}
\end{equation}
We will also write $X_s^{u,t,x}$ instead of $X^u_s$, to
stress dependence on the initial data $t,x$.
By $v:[0,T]\times \calh\to\R$, we denote the mild solution of
the Hamilton-Jacobi-Bellman equation (\ref{HJBformale}).

\begin{theorem}\label{th-rel-font}
 Assume Hypotheses \ref{ipotesiconcrete},
\ref{ipotesiconcretecosto} and \ref{ipotesisupsi}.
 For
every $t\in [0,T]$, $x\in \calh$ and for
 all admissible control $u$ we have $J(t,x,u(\cdot))
 \geq v(t,x)$,
  and the
 equality
$J(t,x,u(\cdot))= v(t,x)$
 holds if and only if
$$
u_s\in \Gamma\left(s,X^{u,t,x}_s, [\nabla
v(\lambda-A)^{1-\beta}](s ,X^{u,t,x}_s)
\;(\lambda-A)^{\beta}D_\lambda\right)
  $$
\end{theorem}

{\bf Proof.} The proof is identical to the proof of relation (7.5)
in \cite[Theorem 7.2]{FuTe}. Just notice that in this case by
(\ref{identifZ}) we have $Z_s^{t,x}= [\nabla
v(\lambda-A)^{1-\beta}](s ,X_s ^{t,x})\;(\lambda-A)^{\beta}D_\lambda$ and
the role of $G$ in \cite[Theorem 7.2]{FuTe} is here played by
$B=(\lambda-A)D_\lambda$. \qed

Under the assumptions of Theorem
\ref{th-rel-font}, let us define the so called
optimal feedback law:
\begin{equation}\label{leggecontrolloottima}
u(t,x)=\gamma\Big(t,x, [\nabla
v(\lambda-A)^{1-\beta}](t ,x)\;(\lambda-A)^{\beta}D_\lambda \Big),\qquad
t\in [0,T],\;x\in \calh.
\end{equation}
Assume that
there exists an adapted process $\{\overline{X}_s,\;s\in
[t,T]\}$ with continuous trajectories solving the
so called closed loop equation: $\P$-a.s.
\begin{equation}\label{cle}
\begin{array}{lll}\dis
\overline{X}_s&=&\dis e^{(s-t)A}x_0+\int_{t}^s
e^{(s-r)A}F(r,\overline{X}_r)\;dr
+
\int_t^s(\lambda -A)^{1-\beta}e^{(s-r)A}(\lambda
-A)^{\beta}D_\lambda\;dW_r
\\&
+&\dis\int_{t}^s(\lambda -A)^{1-\beta}e^{(s-r)A}(\lambda
-A)^{\beta}D_\lambda\;u(r,\overline{X}_r)dr,
\qquad s\in [t,T].
\end{array}
\end{equation}
Then setting
$\overline{u}(s)=u(s,\overline{X}_s)$
we have
$J(t,x,\overline{u}(\cdot))=v(t,x)$ and consequently
 the pair
$(\overline{u},\overline{X})$ is optimal for the control problem.
We nevertheless notice that we do not state conditions
for the existence of a solution of the closed loop equation. Indeed
existence  is not obvious, due to the lack
of regularity of the feedback law $u$ occurring in
(\ref{cle}).

However, under additional assumptions, it is also possible to solve
the closed loop equation (\ref{cle}) and therefore
obtain existence of an optimal control in the present strong
formulation.

We now reformulate the optimal control problem in the
weak sense, following the approach
of \cite{FlSo}. The main advantage is that we
will be able to solve the closed loop
equation, and hence to find an optimal
control, although the feedback law $\underline{u}$
is non-smooth.

\noindent We still assume we are given the
 functions $f$, $\ell$, $\phi$, the corresponding functions
$F$, $\Psi$, $L$, $\Phi$ satisfying Hypotheses \ref{ipotesiconcrete},
\ref{ipotesiconcretecosto} and \ref{ipotesisupsi},
and the set $\calu$
as in the previous sections.
We also assume that initial data $t\in [0,T]$
and $x\in \calh$ are given.
%
%
We call $(\Omega,\calf, ({\cal F}_{t}), \P, W)$
an {\it admissible set-up}, or simply a set-up, if
$(\Omega,\calf, \P)$ is a complete probability space
with a right-continuous and $\P$-complete
filtration $\{{\cal F}_{t},\, t\in [0,T]\}$,
and $\{W_{t},\, t\in [0,T]\}$ is a standard, real valued,
$\calf_t$-Wiener process.

An {\it admissible control
system} (a.c.s.) is defined as
$\U=(\Omega,\calf, ({\cal F}_{t}), \P, W, u, X^u)$ where:
 \begin{itemize}
    \item $(\Omega,\calf, ({\cal F}_{t}), \P, W)$
is an admissible set-up;
    \item $u:\Omega\times [0,T]\to\R$ is an
$({\cal F}_{t})$-predictable process
  with values in $\calu$;
\item $\{X^u_s,\, s\in [t,T]\}$ is an
$({\cal F}_{t})$-adapted continuous process
  with values in $\calh$, mild solution
of the state equation (\ref{eqstatodue})
 with initial condition
 $X^{u}_t=x$.
 \end{itemize}

 By Proposition \ref{regequazforward}, on an arbitrary
 set-up the process $X^u$ is uniquely determined by $u$ and $x$,
 up to indistinguishability.
To every a.c.s.
we associate the cost
$J(t,x,\U)$ given by the right-hand side of
(\ref{costoastratto}).
Although formally the same,
 it is important to note that now the cost
is a functional of the a.c.s., and
not a functional of $u$ alone.
Our purpose is to minimize the functional $J(t,x,\U)$ over all a.c.s.
for fixed initial data $t,x$.
%
%

 \begin{theorem}\label{maincontrollodebole}
 Assume Hypotheses \ref{ipotesiconcrete},
\ref{ipotesiconcretecosto} and \ref{ipotesisupsi}.
 For
every $t\in [0,T]$, $x\in H$,
the infimum of $J(t,x,\U)$ over all a.c.s.
is equal
to $v(t,x)$. Moreover there exists
an a.c.s. $\U=(\Omega,\calf, ({\cal F}_{t}), \P, W, u, X^u)$
for which
 $J(t,x,\U)= v (t,x)$ and the feedback law
$$
u_s= u(s,X^{u}_s) ,\qquad
\P- {\rm a.s. \; for\;a.a.\; } s\in [t,T],
$$
is verified by $u$ and $X^{u}$. Finally,
the optimal trajectory $X^u$ is a weak solution of
the closed loop equation.
 \end{theorem}

{\bf Proof.} We notice that the closed loop equation
(\ref{cle}) always admits a solution in the weak sense 
by an application of the Girsanov theorem.
%
We can apply Theorem \ref{th-rel-font} and obtain all
the required conclusions.
\qed

\section{The forward-backward stochastic differential 
equations in the infinite horizon case}

Eventually we solve the infinite horizon control problem, that we briefly present.
We consider the following infinite horizon cost, with a discount $\mu>0$,
\begin{equation}\label{costoconcretoinfor1}
J(x,u(\cdot))=\E \int_0^{+\infty}e^{-\mu s}\int_0^{+\infty} \ell(\xi,
y(s,\xi),u(s))\;d\xi\;ds .
\end{equation}
that we minimize over all admissible controls. The process
$y$ solves the equation 
\begin{equation}\label{eqconcreta_inf}
  \left\{
  \begin{array}{l}
  \dis
\frac{ \partial y}{\partial s}(s,\xi)= \frac{ \partial^2
y}{\partial \xi^2}(s,\xi)-My(s,\xi)+f(y(s,\xi)), \qquad s\geq 0\;
\xi\in (0,+\infty),
\\\dis
y(0,\xi)=x(\xi),
\\\dis
 y(s,0)= u(s)+\dot{W}_s, \quad
\end{array}
\right.
\end{equation}
where $M>0$ has to be chosen sufficiently large.
Equation (\ref{eqconcreta_inf})
can be reformulated in $\calh$ as
\begin{equation}\label{eqstatoformale_inf}
  \left\{
  \begin{array}{l}
  \dis
dX^u_s= (A-MI)X^u_sds+F(X^u_s)ds+Bu_s ds +B dW_s  \qquad s\geq 0,
\\\dis
X^u_0=x,
\end{array}
\right.
\end{equation}
Also we consider its uncontrolled
version, that is, in mild form,
\begin{equation}\label{eqstatonocontrol_inf}
X_s=\dis e^{s(A-MI)}x+\int_0^s e^{(s-r)(A-MI)}F(X_r^{x})\;dr +
\int_0^se^{(s-r)(A-MI)}B\;dW_r, \qquad s\geq 0.
\end{equation}
By theorem \ref{teosolmild}, for every $T>0$, in $[0,T]$ 
this equation admits a unique mild solution, satisfying
for every $p\in
[1,+\infty)$, $\alpha\in [0,\theta/4)$, 
\begin{equation}\label{stimaes_infor}
    \E\sup_{s\in (0,T]}s^{p\alpha}|X_s^{x}|^p_{D(-A)^\alpha}\leq
c_{p,\alpha}(1+|x|_\calh)^p.
\end{equation}
where $c_{p,\alpha}$ is a constant.
Moreover $X^{x}$ is continuous and G\^ateaux differentiable with respect 
to the initial datum $x$, see proposition \ref{regequazforward}, 
and for every $s\in [0,T]$, we can build the processes 
$\Theta^\alpha(\cdot,x)h$ following proposition \ref{regequazforward2}. 
Moreover $X^x$ admits the Malliavin derivative in every interval $[0,T]$,
see proposition \ref{regequazforwardmalliavin}.
In the next lemma we prove that under our assumptions the derivative $\nabla_xX^x_\cdot$
and the process $\Theta^\alpha(\cdot,x)h$ are uniformly bounded in time.
\begin{lemma}\label{lemmader_forward_lim}
Assume that hypothesis \ref{ipotesiconcrete} holds true and that in
equation (\ref{eqstatonocontrol_inf}) $M$ is sufficiently large (to be chosen in the following proof), then there exists 
a constant $C>0$ such that 
$\vert\nabla_xX^x_t\vert+ \vert\Theta^\alpha(t,x)h\vert\leq C\vert h\vert $
for every  $t>0$ and every $x,h\in\calh$.
\end{lemma}
{\bf Proof.}
We already know that, see proposition \ref{regequazforward},
the map $x\to X_\cdot^{x}$
has, at every point $x\in \calh$, in every direction $h\in\calh$, 
a G\^{a}teaux derivative $\nabla_xX_\cdot^{x}h$ 
and the map $x\to \nabla_x X_\cdot^x$ belongs to
$\calg^1(\calh, L^p_\calp(\Omega;C([0,T];\calh))$ 
and, for every direction $h\in \calh$, the following
equation holds $\P$-a.s.:
\begin{equation}\label{derdixmild_inf}
    \nabla_x X_t^{x}h=e^{t(A-MI)}h+
    \int_0^t e^{(t-s)(A-MI)}\nabla_xF(
    X_s^{x}) \nabla_x X_s^x\; ds
    \qquad t\geq 0.
\end{equation}
Since by \cite{Kr01}, theorem 2.5, there exists $C>0$, independent
on $f$, such that for every $f\in \calh$
$$
\vert e^{tA}f \vert_\calh \leq C \vert f \vert_\calh, \qquad t\geq 0,
$$
then 
$$
\vert e^{t(A-MI)}f \vert_\calh \leq C e^{-Mt}\vert f \vert_\calh, \qquad t\geq 0.
$$
So, by equation (\ref{derdixmild_inf}), we can deduce that 
\begin{align*}
\vert\nabla_x X_t^{x}h\vert
&\leq \vert e^{t(A-MI)}h\vert+ \vert\int_0^t e^{(t-s)(A-MI)}\nabla_xF(
    X_s^{x}) \nabla_x X_s^x\; ds\vert \\ \nonumber
 &\leq Ce^{-Mt}\vert h\vert+C_f C\sup_{0\leq s\leq t}\vert\nabla_x X_s^{x}h\vert 
\int_0^t e^{-M(t-s)}ds \\ \nonumber
&  \leq C\left[ \vert h\vert+\frac{C_f}{M} \sup_{0\leq s\leq t}\vert  X_s^{x} \vert\right]\leq  
  C\left[ \vert h\vert+\frac{C_f}{M} \sup_{s\geq 0}\vert \nabla_x X_s^{x} \vert\right]  \\ \nonumber
\end{align*}
The previous inequality holds true for every $t>0$, so we get
\[
   \sup_{t\geq 0} \vert \nabla_x X_t^{x}h\vert \leq C\left[
 \vert h\vert+\frac{C_f}{M} \sup_{t\geq 0}\vert  X_t^{x} \vert\right], 
\]
and assuming that $M>C_f \times C$ we obtain that
$\vert\nabla_x X_s^x h\vert\leq C \vert h\vert$
for some constant $C>0$ and for every $x,h\in \calh$. 

\noindent For what concerns $\Theta^{\alpha}$, we already know that they satisfy 
and equation like \ref{equazpertheta}, with $A-MI$ in the place of $A$. 
We also remark that for every $\gamma>0$,
$\label{stima_su_B}
\vert e^{(tA)B} \vert\leq M_\gamma t^{-\gamma}$, so it follows that
\begin{align*}
\vert\Theta^{\alpha}(t,x)h\vert
& \leq C\times C_f 
\sup_{0\leq s\leq t}\vert \Theta^{\alpha}(s,x)h \vert\int_0^t e^{-M(t-s)}ds 
+ M_\alpha C\times C_f \vert h \vert\int_0^t e^{-M(t-s)} s^{-\alpha}ds\\ \nonumber
& \leq C\times C_f \left[
\frac{1}{M}\sup_{0\leq s\leq t}\vert \Theta^{\alpha}(s,x)h  
+ \frac{1}{1-\alpha}M_\alpha\vert h \vert +\frac{M_\alpha}{M}\vert h \vert\chi_{[1,+\infty)}(t)\right].\\ \nonumber
\end{align*}
As for $\nabla_x X$, this inequality holds true for every $t>0$ and again
if $M>C_f \times C$ we obtain that
$\vert\nabla_x X_s^x h\leq C \vert h\vert$
for some constant $C>0$ and for every $x,h\in \calh$
\qed

From now on, and in equations (\ref{eqconcreta_inf}), (\ref{eqstatoformale_inf}) and 
(\ref{eqstatonocontrol_inf}), we take $M> C\times C_f$.

We need to notice that
proposition \ref{ucompostox} can be easily adequated to
the case of a function $w$ not depending on time, so
we can state the following result about the joint quadratic variation of
the process $u(X_\cdot)$ with $W$, where $u is a suitable function$
\begin{proposition}\label{ucompostox_rem} 
Suppose that $w\in C( \calh;
\mathbb{R})$ is G\^{a}teaux differentiable 
and that there exist constants $K$ and $m$ such that
\begin{equation*}
 |w(x)|\leq K(1+|x|)^m,
\qquad |\nabla w (x) |\leq K(1+|x|)^m,\qquad
x\in H.
\end{equation*}
Assume that for every
$x\in \calh$, $\beta \in (0,\frac{1}{2}+\frac{\theta}{4})$, the linear
operator $k\to \nabla w(x)(\lambda-A)^{1-\beta}k$ (a priori
defined for $k\in D(\lambda-A)^{1-\beta}$)
has an extension to a bounded linear
operator $\calh\to \mathbb{R}$, that we
 denote by $[\nabla w(\lambda -A)^{1-\beta}](x)$.
Moreover assume that
the map $(x,k)\to [\nabla w(\lambda-A)^{1-\beta}](x)k$
is continuous $ \calh\times \calh \to \mathbb{R}$.

\noindent For $x\in \calh$, let
$\{X_s^{x},\,s>0\}$ be the solution of
equation (\ref{eqstatonocontrol_inf}).
Then the process $\{w(X_s^{x}),\, s>0\}$
admits a joint quadratic variation process
 with $W$, on every interval
$[t,s]\subset [0,+\infty)$, given by
  $$
  \int_t^s[\nabla w(\lambda -A)^{1-\beta}](X_r^{x})\,
(\lambda -A)^{\beta}D_\lambda\;dr.
$$
\end{proposition}

\noindent In order to solve the infinite horizon control problem, we 
consider the following backward stochastic
differential equation:
\begin{equation}\label{BSDE_inf}
  dY^x_s=-\Psi(X_s^x,Z^x_s)\;ds+\mu Y^x_s\;ds+Z^x_s \;dW_s,\qquad s\geq 0,\\
\end{equation}
for the unknown real processes $Y^x$, $Z^x$, also denoted by $Y$ a nd $Z$.
The equation is understood in the usual way: $\P$-a.s., for
every $T>0$,
\begin{equation}\label{BSDEmild_inf}
Y^x_s+\int_s^TZ^x_r \;dW_r=
Y^x_T+\int_s^T\left(\Psi(X_r^x,Z^x_r)-\mu Y^x_r\right)dr,\qquad s\geq 0.
\end{equation}
We make the following assumptions:
\begin{hypothesis}\label{ipsupsi_inf}
\begin{itemize}
 \item [i)]The function $\Psi : \calh \times \R \to \R$ is continuous
in $x$ and uniformly Lipschitz continuous in $z$ that is 
$\vert \Psi(x,z_1)-\Psi(x,z_2)\vert\leq K \vert z_1-z_2 \vert $
\item [ii)] $\sup_{x\in\calh}\vert \Psi(x,0) \vert :=M<+\infty $
\item [iii)] $\mu>0$.
\end{itemize}
\end{hypothesis}

We can state the following result on existence and uniqueness of
a solution $(Y,Z)$ of equation (\ref{BSDE_inf}).

\begin{proposition}\label{propBSDE_inf}Assume hypotheses 
\ref{ipotesiconcrete} and \ref{ipsupsi_inf}, 
\begin{itemize}
 \item  [i)]For any $x\in\calh$ equation (\ref{BSDE_inf})
admits a unique solution $(Y^x,Z^x)$ such that $Y^x$ is a 
continuous process bounded by $M/\mu$, and $Z\in L^2_{\calp,loc}((0,+\infty),\R)$
with $\E\int_0^{+\infty}e^{-2\mu s}\vert Z_s\vert^2ds<\infty$.
The solution is unique in the class of processes such that $Y^x$
is continuous and bounded and $Z^x\in L^2_{\calp,loc}((0,+\infty),\R)$.
\item [ii)] Denoting by $(Y^{n,x},Z^{n,x})$ the solution to the following, finite horizon, 
BSDE
\begin{equation}\label{BSDE_inf_approx}
 Y_s^{n,x}+\int_s^n Z_r^{n,x} \;dW_r=
+\int_s^n\left(\Psi(X_r^x,Z_r^{n,x})-\mu Y_r^{n,x}\right)dr,
\end{equation}
then $\vert Y_s^{n,x}\vert\leq \frac{M}{\mu}$ and the following convergence rate holds:
\begin{equation}\label{BSDE_convergencerate}
 \vert Y_s^{n,x}-Y_s^{x} \vert\leq \frac{M}{\mu}e^{-\mu(n-s)}.
\end{equation}
Moreover
\begin{equation}\label{BSDE_convergencerateZ}
\E\int_0^{+\infty}e^{-2\lambda s} \vert Z_s^{n,x}-Z_s^{x}\vert ds \to 0.
\end{equation}
\item [iii)] For all $T>0$ and $p\geq 1$, the map $x\to
(Y^{x}\lvert_{[0,T]},Z^{x}\lvert_{[0,T]})$ is continuous from $\calh$
to $L^p_{\calp}(\Omega,C([0,T],\R)\times L^p_{\calp}(\Omega,L^2([0,T],\R)$
\end{itemize}
\end{proposition}

{\bf Proof.}
The proof follows the proof of proposition 3.2 in \cite{HuTess}, in the case of $\psi$
not depending on $Y$. 
\qed

\noindent We assume moreover the following:
\begin{hypothesis}\label{ip_psiagg_inf}
$\Psi \in \calg^1(\calh\times \R)$ and $\nabla_x\Psi(x,z)\leq c$ for
every $x\in\calh$, $z\in\R$, and for some constant $c>0$. $\big(\nabla_z\Psi(x,z)\leq c$
is also bounded as a consequence of hypothesis \ref{ipsupsi_inf}, point i) $\big)$.
\end{hypothesis}
We can state the following theorem:
\begin{theorem}\label{bsderegolare-inf}
Assume that hypotheses \ref{ipotesiconcrete}, \ref{ipsupsi_inf} and \ref{ip_psiagg_inf} hold true.
Then the map $x\to Y_0^x \in \calg ^1(\calh, \R)$ and $\vert Y_0^x +\nabla Y_0^x\vert \leq C$.
\end{theorem}
{\bf Proof.}
We follow the proof of Theorem 3.1 in \cite{HuTess}. In that theorem,
it was assumed that the operator $A+\nabla_xF(x)$ is dissipative.
This is used in order to prove, see Lemma 3.1 in \cite{HuTess},
that $\vert\nabla_x X_s^x h\leq C \vert h\vert$
for some constant $C>0$ and for every $x,h\in \calh$. In the present situation,
we already know that, see lemma \ref{lemmader_forward_lim},
$\vert\nabla_x X_s^x h\leq C \vert h\vert$
for some constant $C>0$ and for every $x,h\in \calh$.


\noindent The proof now follows exactly from the proof of theorem 3.1 in 
\cite{HuTess}.
\qed

Next we have to prove a further regularity result, similar to the one stated in proposition \ref{regback}.
To this aim, we need to adapt the results in \cite{BriHu}.
\begin{lemma}\label{lemma2.1}
Let us consider the following BSDE on an infinite horizon,
\begin{equation}\label{BSDEHuBri}
 Y_t=Y_T-\int_t^T Z_r dW_r+\int_t^T\left(f_1(r,Z_r)+f_2(r)-\mu Y_r\right)dr,
\end{equation}
where $\mu>0$, $f_1:\Omega \times [0,+\infty) \times \R\rightarrow\R $ and 
$\forall z\in\R$ the process $(f_1(t,z))_{t\geq0}$ is predictable, and the process $(f_2(t))_{t\geq0}$ is predictable. Moreover assume that:
\begin{itemize}
 \item  [i)]$f_1$ is uniformly lipschit continuous in $z$ with lipschitz 
constant $K$: $\forall\, t>0$, $\forall\, z_1,z_2\in\R$ 
\[
 \vert f_1(t,z_1)-f_1(t,z_2) \vert\leq K\vert z_1-z_2 \vert, \qquad \P-a.s.,
\]
and $f_1(t,0)$ is bounded.
\item [ii)] there exists a constant $M>0$ and a function $g\in L^1([0,1], \R)$, such that for every $t\geq 0$,
\[
  \vert f_2(t) \vert\leq \vert g(t)\vert\chi_{[0,1]}(t)+M.
\]
\end{itemize}
Then 
\begin{itemize}
\item  [i)]there exists a solution $(Y,Z)$ to equation (\ref{BSDEHuBri}) such that $Y$ is a continuous, predictable
process, bounded by a constant $C>0$ and $Z\in L^2_{\calp, loc}([0,+\infty),\R)$ and the solution is unique 
in such class of processes. Moreover
$\int_0^{+\infty}e^{-2\mu s}\vert Z_s \vert^2 <+\infty$.
\item [ii)]Denoting by $(Y^n, Z^n)$ the unique solution of the BSDE
\begin{equation} \label{BSDEHuBri_n}
 Y^n_t=-\int_t^n Z^n_r dW_r+\int_t^n\left(f_1(r,Z^n_r)+f_2(r)-\mu Y^n_r\right)dr,
\end{equation}
\end{itemize}
then $\vert Y^n_t\vert \leq C$ and the following convergence rate holds:
\[
 \vert Y^n_t-Y_t \vert \leq C e^{-\mu(n-t)}
\]
and moreover
\[
 \int_0^{+\infty} e^{-2\mu t}\vert Z^n_t-Z_t \vert^2 dt \to 0.
\]
\end{lemma}
{\bf Proof.}
Let us consider $(Y^n, Z^n)$ solution to equation (\ref{BSDEHuBri_n}). 
We set 
\[
 f_3(r):=
\left\lbrace
\begin{array}[l]{ll}
 \dfrac{f_1(r,Z^n_r)-f_1(r,0)}{\vert Z^n_r \vert^2}Z^n_r & \text{ if }Z^n_r\neq0 \\
0 & \text{ otherwise }.
\end{array}
\right. 
\]
By our assumptions on $f_1$, $f_3$ is bounded and so by the Girsanov theorem,
there exists a probabilty measure $\tilde{\P}$, equivalent to the original one $\P$,
such that 
\[
 \left\lbrace \tilde{W}_t=-\int_0^t f_3(r) dr+W_t, t\geq 0 \right\rbrace 
\]
is a Brownian motion.
So in $(\Omega, \calf, \tilde{\P})$ equation (\ref{BSDEHuBri_n}) can be rewritten as
\begin{equation*} 
 Y^n_t=-\int_t^n Z^n_r d\tilde{W}_r+\int_t^n\left(f_1(r,0)+f_2(r)-\mu Y^n_r\right)dr,
\end{equation*}
Since $f_1(\cdot,0)$ is bounded and $f_2$ is 
integrable near $0$ and bounded 
otherwise, by the Gronwall lemma it 
follows that for every $t\in[0,n]$
\[
 \vert Y^n_t \vert \leq C e^{-\mu (n-t)}\tilde{\E}^{\calf _t}
\int _t^n (\vert f_1(r,0)\vert +\vert f_2(r)\vert dr
\]
where $C$ is a constant independent on $n$.
By setting $Y$ as the pointwise limit 
of $Y^n$ we get that $Y$ is bounded.
By applying It\^o formula to $e^{-2\mu t}\vert Y^n_t\vert ^2$
it follows that
\[
 \int_0^{+\infty} e^{-2\mu t}\vert Z^n_t \vert^2 dt < + \infty.
\]
Now let us define $\tilde{Y}_t^n:=Y^n_t-Y_t$ and $\tilde{Z}_t^n:=Z^n_t-Z_t$. $(\tilde{Y}_t^n,\tilde{Z}_t^n)$
solve, for $t\in [0,n]$, the following BSDE:
\begin{equation}\label{tildeBSDE}
 \tilde{Y}^n_t=-Y_n-\int_t^n \tilde{Z}^n_r dW_r+\int_t^n\left(f_1(r,Z^n_r)-f_1(r,Z_r)\right)dr-
\int_t^n\mu\tilde{Y}^n_rdr 
\end{equation}
We also set 
\[
 F_r:=
\left\lbrace
\begin{array}[l]{ll}
 \dfrac{f_1(r,Z^n_r)-f_1(r,Z_r)}{\vert Z^n_r-Z_r \vert^2}(Z^n_r-Z_r) & \text{ if }Z^n_r-Z_r\neq0 \\
0 & \text{ otherwise }.
\end{array}
\right. 
\]
By the lipschitz assumptions on $f_1$, $F$ is bounded and so by the Girsanov theorem,
there exists a probabilty measure $\bar\P$, equivalent to the original one $\P$,
such that 
\[
 \left\lbrace \bar W_t=-\int_0^t F_r dr+W_t, t\geq 0 \right\rbrace 
\]
is a Brownian motion.
In $(\Omega, \calf, \bar P)$ equation (\ref{tildeBSDE}) can be rewritten as
\begin{equation*}
 \bar Y^n_t=-Y_n-\int_t^n \tilde{Z}^n_r d\tilde{W}_r-
\int_t^n\mu\tilde{Y}^n_rdr 
\end{equation*}
So the following rate of convergence holds true:
\[
\vert Y_t-Y^n_t \vert \leq C e^{-\mu (n-t)},
\]
where $C$ is a constant not depending on $n$.
By applying It\^o formula to $e^{-2\mu t}\vert\tilde{Y}^n_t\vert^2$
it follows that 
\[
 \int_0^{+\infty} e^{-2\mu t}\vert Z^n_t-Z_t \vert^2 dt \to 0.
\]
\qed

\begin{theorem}\label{regback_inf}
For every $\alpha \in [0,1/2)$, $p\in [2,\infty)$ there exist two
functions $P^{\alpha}(x)k$ and $ Q^{\alpha}(x)k$, $x\in \calh$, $k\in \calh$
such that if
$k\in D(\lambda-A)^{\alpha}$, $x\in \calh$, then
\begin{equation}\label{def-di-P-inf}
    P^{\alpha}(x)k=  \nabla_x Y_0^x(\lambda-A)^{\alpha}k \\
 \end{equation}
and
\begin{equation}\label{def-di-Q-inf}
Q^{\alpha}(x)k= \nabla_x Z_0^{x}(\lambda-A)^{\alpha}k 
\end{equation}
Moreover the map $ (x,k)\to P^{\alpha}(x)k$ 
is continuous from $\calh$ to
$\mathbb{R}$ and 
linear with respect to $k$.

Finally there exists a constant $C_{\nabla Y,\alpha,p}$ such that
\begin{equation}\label{stimadiP_inf}
 |P^{\alpha}(x)k| \leq
C_{\nabla Y,\alpha} |k|_\calh.
\end{equation}
\end{theorem}
{\bf Proof.}
For $x\in \calh$ and $k\in D(\lambda-A)^{\alpha}$, let
$P^{\alpha}(x)k$ and $Q^{\alpha}(x)k $ be defined by
(\ref{def-di-P-inf}) and (\ref{def-di-Q-inf}) respectively.
By Theorem \ref{bsderegolare-inf}
the map $k\to(P^{\alpha}(x)k,Q^{\alpha}(x)k )$ is
a bounded linear operator from $D(\lambda-A)^{\alpha}$ to
$\mathbb{R}\times\mathbb{R}$. 

\noindent Let us introduce the pair of processes
$(P^{\alpha}(\cdot, x)k,Q^{\alpha}(\cdot,x)k)$
solution of the following BSDE
\begin{align}\label{BSDE_PQ^alpha_inf}
 P^{\alpha}(t,x)k& =P^{\alpha}(T,x)k+
\int_t^T \nabla_x\Psi(X_s^x,Z_s^{x})
\left( \Theta^\alpha(s,x)k+(\lambda-A)^\alpha e^{sA} k \right)ds\\ \nonumber
& -\int_0^t \mu P^{\alpha}(s,x)k\;ds 
+\int_t^T \left[\nabla_z\Psi(X_s^y,Z_s^{x})Q^{\alpha}(s,x)k\right]ds 
- \int_t^T Q^{\alpha}(s,x)k\; dW_s,\qquad t\geq 0. 
\end{align}
Equation (\ref{BSDE_PQ^alpha_inf}) admits a unique bounded solution 
by applying lemma \ref{lemma2.1}.

\noindent Moreover, let us define the processes
$(P^{\alpha,n}(\cdot,x)k,Q^{\alpha,n}(\cdot,x)k $ solution of the equation
\begin{align}\label{equaz_Pn_inf}
  P^{\alpha,n}(t,x)k&=\int_t^n \nabla_x\Psi(X_s^x,Z_s^{n,x})
  \left( \Theta^\alpha(s,x)k+(\lambda-A)^\alpha e^{sA} k \right)\;ds 
-\int_t^n \mu P^{\alpha,n}(s,x)\;ds \\ \nonumber
  &+\int_t^n\nabla_z \Psi(X_s^x,Z^{n,x}_s)\;Q^{\alpha,n}(s,x)k\; ds
 + \int_t^nQ^{\alpha,n}(s,x)k\; dW_s,\qquad t\geq 0 
\end{align}
We notice that equation (\ref{equaz_Pn_inf}) is obtained by
formally deriving equation (\ref{BSDE_inf_approx}) in the direction $(\lambda-A)^\alpha k$.

\noindent Equation (\ref{equaz_Pn_inf}) can be rewritten
\begin{align}\label{equaz_Pn_inf_2}
 P^{\alpha,n}(t,x)k&=\int_t^n \nu^n(s,x)k\;ds -\int_t^n \mu P^{\alpha,n}(s,x)k\;ds\\ \nonumber
  &+\int_t^n\nabla_z \Psi(X_s^x,Z_s^{n,x})\;Q^{\alpha,n}(s,x)k\; ds
 + \int_t^nQ^{\alpha,n}(s,x)k\; dW_s,\qquad 0\leq t\leq n,
\end{align}
where
$$
  \nu^n(s,x)k=   \nabla_x\Psi(X_s^x,Z_s^{n,x})
  \left(\Theta^{\alpha} (s,x)k+ (\lambda-A)^{\alpha} e^{sA}k\right),
\qquad t\in[0,n].
$$
Now we choose arbitrary $k\in \calh$ and notice that $\nu^n(s,x)k$
can still be defined by the above formulae. Hypothesis
\ref{ip_psiagg_inf} and relation (\ref{stimadi nabla theta}) yield:
\begin{align*}
   \E\left(\int_0^n |\nu^n (s,x)k|^2ds \right)^{p/2}&\leq c_1
\E\!\left(\int_0^n \left( |\Theta^{\alpha}
(s,x)k| +s^{-\alpha}|k|_\calh\right)^2ds\right)^{p/2} \\ \nonumber
  &
  \dis c_2 \left[n^{p/2}+n^{(1-2\alpha)p/2}\right] |k|^p
  \leq c_3|k|^p,
\end{align*}
where $c_1$, $c_2$ and $c_3$ are suitable constants independent on
$t,x,k$.
By Proposition 4.3 in \cite{FuTe}, for all $k\in \calh$ there
exists a unique pair
$((P^{\alpha,n}(\cdot,x)k,Q^{\alpha,n}(\cdot,x)k )$
 belonging to
$L^p_{\mathcal{P}}(\Omega,C([0,n],\mathbb{R}))\times
L^p_{\mathcal{P}}(\Omega,L^2([0,n],\mathbb{R}))$ and solving
equation (\ref{equaz_Pn_inf_2}).
Moreover, by applying proposition \ref{regback}, we get that 
$ (x,k)\to P^{\alpha,n}(\cdot,x)k$ 
and the map $(x,k)\to Q^{\alpha,n}(\cdot,x)k $
are continuous from $[0,n)\times \calh\times \calh$ to
$L^p_{\mathcal{P}}(\Omega,C([0,n],\mathbb{R}))$ and 
linear with respect to $k$.
Finally there exists a constant $C_{\nabla Y,\alpha,p,n}$ such that
\begin{equation}\label{stimadiP^n}
    \E\sup_{s\in [0,n]}|P^{\alpha,n}(s,x)k|_\calh^p+\E\left(\int_0^n
|Q^{\alpha,n} (s,x)k|_{(\mathbb{R})} ds\right)^{p/2} \leq
C_{\nabla Y,\alpha,p,n} n^{-\alpha p} (1+|x|_\calh)^p|k|_\calh^p.
\end{equation}
Moreover we want to prove that $P^{\alpha,n}(\cdot,x)k$ is 
a bounded process, uniformly in $n$.
To this aim let $k\in D(\lambda-A)^\alpha$ and let $x, y\in \calh $ such that
$x-y=(\lambda-A)^\alpha k $.

\noindent Let us also define $\tilde{Y}^{n}_t=Y^{n,x}_t-Y^{n,y}_t $ and 
$\tilde{Z}^{n}_t=Z^{n,x}_t-Z^{n,y}_t$. So the pair $(\tilde{Y}^{n},\tilde{Z}^{n})$ solves 
the following backward stochastic differential equation:
\begin{equation*}
  \tilde{Y}^n_t=\int_t^n \left[\Psi(X_s^x,Z_s^{n,x})-\Psi(X_s^y,Z_s^{n,y})\right]ds
   -\int_t^n \mu \tilde{Y}^n_s\;ds
 + \int_t^n\tilde{Z}^{n}_s\; dW_s,\qquad t\in [0,n],
\end{equation*}
that we can also write as 
\begin{align*}
\tilde{Y}^n_t & =\int_t^n \left[\Psi(X_s^x,Z_s^{n,x})-\Psi(X_s^y,Z_s^{n,x})\right]ds
    -\int_t^n \mu \tilde{Y}^n_s\;ds \\ \nonumber
& +\int_t^n \left[\frac{\Psi(X_s^y,Z_s^{n,x})-\Psi(X_s^y,Z_s^{n,y})}
{\vert Z_s^{n,x}-Z_s^{n,y} \vert^2} \left( Z_s^{n,x}-Z_s^{n,y}\right)\right]\tilde{Z}^{n}_sds 
 + \int_t^n\tilde{Z}^{n}_s\; dW_s,\qquad t\in [0,n].\\ \nonumber
\end{align*}
Since $\Psi$ is uniformly lipschitz with respect to $z$, 
\[
\left\vert\frac{\Psi(X_s^y,Z_s^{n,x})-\Psi(X_s^y,Z_s^{n,y})}
{\vert Z_s^{n,x}-Z_s^{n,y} \vert^2}\left( Z_s^{n,x}-Z_s^{n,y}\right)\right\vert
\leq L. 
\]
So, as in \cite{BriHu}, lemma 3.1, by the 
Girsanov theorem there exists a probability measure $\tilde{\P}$ such that in 
$(\Omega, \calf, \tilde{\P})$ the process $(\tilde{W}_t)_{t\in[0,n]}$ defined by
\[
\tilde{W}_t= W_t+
\int_0^t  \frac{\Psi(X_s^y,Z_s^{n,x})-\Psi(X_s^y,Z_s^{n,y})}
{\vert Z_s^{n,x}-Z_s^{n,y} \vert^2}\left(Z_s^{n,x}-Z_s^{n,y}\right)  ds 
\]
is a real brownian motion. In this probability space, $(\tilde{Y}^{n},\tilde{Z}^{n})$ solve 
the following backward stochastic differential equation:
\begin{equation*}
  \tilde{Y}^n_t=\int_t^n \left[\Psi(X_s^x,Z_s^{n,x})-\Psi(X_s^y,Z_s^{n,x})\right]ds
   -\int_t^n \mu \tilde{Y}^n_s\;ds
 + \int_t^n\tilde{Z}^{n}_s\; d\tilde{W}_s,\qquad t\in [0,n],
\end{equation*}
Taking the conditional expectation in the previous equation,
we get that 
\[
 \tilde{Y}^n_t=\tilde{\E}^{\calf_t}
\int_t^n \vert\Psi(X_s^x,Z_s^{n,x})-\Psi(X_s^y,Z_s^{n,x})\vert ds
-\mu\tilde{\E}^{\calf_t}\int_t^n \tilde{Y}^n_s ds.
\]
By lemma \ref{lemmader_forward_lim}, and since $\Psi$ is lipschitz 
with respect to $x$, we get
\begin{align*}
 \vert\tilde{Y}^n_t\vert & \leq \tilde{\E}^{\calf_t}
\int_t^ne^{\mu(t-s)} \vert\Psi(X_s^x,Z_s^{n,x})-\Psi(X_s^y,Z_s^{n,x})\vert ds \\ \nonumber
&\leq L\int _t^ne^{\mu(t-s)} \vert e^{sA}(\lambda-A)^\alpha k+\Theta^\alpha(s,x)k\vert ds \\ \nonumber
&\leq C L\int _t^ne^{\mu(t-s)} \left(s^{-\alpha}+1\right) \vert k\vert_\calh ds \leq C\vert k\vert_\calh,
\end{align*}
where $C$ is a constant that may change its value from line to line, and that 
does not depend on $n$. So we have
\[
\sup_{t\in[0,n]}\vert\tilde{Y}^n_t\vert  \leq C \vert k\vert_\calh,
\]
and consequently we get that for every $x,k\in \calh$,
\[
 \sup_{t\in[0,n]}\vert P^{\alpha,n}(t,x)k\vert \leq C \vert k\vert_\calh.
\]
Then, again as in \cite{BriHu}, by applying the It\^o formula to 
$e^{-2\mu t}\vert P^{\alpha,n}(t,x)k\vert^2$, we get
\begin{align*}
 d e^{-2\mu t}\vert P^{\alpha,n}(t,x)k\vert^2& =
-2\mu e^{-2\mu t}\vert P^{\alpha,n}(t,x)k\vert^2dt+
2\mu e^{-2\mu t}\vert P^{\alpha,n}(t,x)k \vert^2 dt \\ \nonumber
&-2e^{-2\mu t} \nu^n(t,x)kP^{\alpha,n}(t,x)kdt
+ e^{-2\mu t} \vert Q^{\alpha,n}(t,x)k \vert^2 dt
 \\ \nonumber
&-2e^{-2\mu t}Q^{\alpha,n}(t,x)kP^{\alpha,n}(t,x)kdW_t 
-2e^{-2\mu t}\nabla_z\Psi(X_t^x,Z_t^{n,x})Q^{\alpha,n}(t,x)k P^{\alpha,n}(t,x)k dt .\nonumber
\end{align*}
By taking expectation, we get 
\begin{align*}
\E e^{-2\mu t}\vert P^{\alpha,n}(t,x)k\vert^2& =
\E\int _t^n e^{-2\mu s} \nu^n(s,x)P^{\alpha,n}(s,x)kds
-\int_t^n e^{-2\mu s} \vert Q^{\alpha,n}(s,x)k \vert^2ds \\ \nonumber
&+2\int_t^n e^{-2\mu s}\nabla_z\Psi(X_s^x,Z_s^{n,x})Q^{\alpha,n}(s,x)k P^{\alpha,n}(s,x)k ds 
\end{align*}
By Young inequality and since $P^{\alpha,n}(\,x)k$ is a uniformly bounded process we get that
\[
 \E\int_0^{+\infty}e^{-2\mu t}\left( \vert P^{\alpha,n}(t,x)k\vert^2 +\vert Q^{\alpha,n}(t,x)k\vert^2\right)dt<+\infty.
\]
Now our proof substantially follows the proof of theorem 3.1 in \cite{HuTess}.
Let $\calm ^{2, -2\mu}$ be the Hilbert space of all couples of real valued, $(\calf_t)_{t\geq0}$-adapted 
processes $(y,z)$, such that
\[
 \vert (y,z)\vert_{\calm ^{2, -2\mu}}^2=
\E\int_0^{+\infty}e^{-2\mu t}
\left( \vert y_t\vert^2 +\vert z_t\vert^2\right)dt<+\infty.
\]
Fixed $x,k\in\calh$, there exists a subsequence of 
$(P^{\alpha,n}(\cdot,x)k,Q^{\alpha,n}(\cdot,x)k,P^{\alpha,n}(0,x)k)$, which we still denote 
by itself, such that $(P^{\alpha,n}(\cdot,x)k,Q^{\alpha,n}(\cdot,x)k,P^{\alpha,n}(0,x)k)$
converges weakly in $\calm ^{2, -2\mu}$ 
to $(U^{1,\alpha}(\cdot,x)k,V^{1,\alpha}(\cdot,x)k,\xi (x,k))$.

\noindent Next we define 
\begin{align}\label{eqU2V2}
  U^{2,\alpha}(t,x)k&=\xi(x,k)-\int_0^t \nabla_x\Psi(X_s^x,Z_s^{x})
\left( \Theta^\alpha(s,x)k+e^{sA}(\lambda-A)^\alpha k \right)ds
 -\int_0^t \mu U^{1,\alpha}(s,x)k\;ds \\ \nonumber
&-\int_0^t \left[\nabla_z\Psi(X_s^y,Z_s^{x})V^{1,\alpha}(s,x)k\right]ds 
+ \int_0^t V^{1,\alpha}(s,x)k\; dW_s,\qquad t\geq 0
\end{align}
where $(Y,Z)$ is the unique bounded solution of equation (\ref{BSDE_inf}).
Let us rewrite, for $t\in[0,n]$, equation (\ref{equaz_Pn_inf}) as
\begin{align*}
P^{\alpha,n}(t,x)k & =P^{\alpha,n}(0,x)k-\int_t^n \nabla_x\Psi(X_s^x,Z_s^{n,x})
  \left( \Theta^\alpha(s,x)k+e^{sA}(\lambda-A)^\alpha k \right)\;ds \\ \nonumber
&+\int_0^t \mu P^{\alpha,n}(s,x)\;ds
 +\int_0^t\nabla_z \Psi(X_s^x,Z^{n,x}_s)\;Q^{\alpha,n}(s,x)k\; ds
 + \int_t^nQ^{\alpha,n}(s,x)k\; dW_s. \nonumber
\end{align*}
As in \cite{HuTess}, theorem 3.1, we can deduce that $P^{\alpha,n}(\cdot,x)k$ coverges
weakly to $U^{2,\alpha}(\cdot,x)k$ in $L^2_\calp([0,T];\R)$.
Moreover, by lemma \ref{lemma2.1}, $(U^{2,\alpha}(\cdot,x)k,V^{1,\alpha}(\cdot,x)k)$ is the unique bounded solution
to equation 
\begin{align*}
  U^{2,\alpha}(t,x)k&=U^{2,\alpha}(0,x)k-\int_0^t \nabla_x\Psi(X_s^x,Z_s^{x})
\left( \Theta^\alpha(s,x)k+e^{sA}(\lambda-A)^\alpha k \right)ds
 +\int_0^t \mu U^{2,\alpha}(s,x)k\;ds \\ \nonumber
&-\int_0^t \left[\nabla_z\Psi(X_s^y,Z_s^{x})V^{1,\alpha}(s,x)k\right]ds 
+ \int_0^t V^{1,\alpha}(s,x)k\; dW_s,\qquad t\geq 0,\\ \nonumber
\end{align*}
so we also have $(U^{2,\alpha}(\cdot,x)k,V^{1,\alpha}(\cdot,x)k)=
(P^{\alpha}(\cdot,x)k,Q^{\alpha}(\cdot,x)k)$, where $(P^\alpha,Q^\alpha)$
solve BSDE (\ref{BSDE_PQ^alpha_inf}),
and in particular $U^{2,\alpha}(0,x)k=\xi(x)k$ is the limit of
$P^{\alpha,n}(0,x)k$ along the original sequence. 

Now we have to prove that the map $x\to U^{2,\alpha}(0,x)k$ is continuous.
Let us consider $(U^{n,\alpha}, V^{n,\alpha})$ the unique solution of equation
\begin{align*}
  U^{n,\alpha}(t,x)k&=\int_t^n \nabla_x\Psi(X_s^x,Z_s^{x})
\left( \Theta^\alpha(s,x)k+e^{sA}(\lambda-A)^\alpha k \right)ds
 -\int_t^n \mu U^{n,\alpha}(s,x)k\;ds \\ \nonumber
&+\int_t^n \left[\nabla_z\Psi(X_s^y,Z_s^{x})V^{n,\alpha}(s,x)k\right]ds 
-\int_t^n V^{n,\alpha}(s,x)k\; dW_s,\qquad t\geq 0. \nonumber
\end{align*}
By proposition \ref{regback}, the map $x\to U^{n,\alpha}(0,x)k$ is 
continuous, and by arguments similar to the ones used before $U^{n,\alpha}(\cdot,x)k$
is a uniformly bounded process. In the probability space 
$(\Omega, \calf, \hat{\P})$ where $\hat{\P}$ is a probability measure,
equivalent to $\P$, such that the process 
\[
\left\lbrace \hat{W}_t:= -\int_0^t \nabla_z\psi(X_s^y,Z_s^{x})ds 
+W_t,\qquad t\geq 0. \right\rbrace
\]
is a Brownian motion, $(e^{-\mu t}U^{n,\alpha}(t,x)k-U^{2,\alpha}(t,x)k),
e^{-\mu t}(V^{n,\alpha}(t,x)k-V^{1,\alpha}(t,x)k)_{t\in[0,n]}$
solve the following BSDE,
\[
\left\lbrace
\begin{array}{l}
  de^{-\mu t}(U^{n,\alpha}(t,x)k-U^{2,\alpha}(t,x)k )=
e^{-\mu t}(V^{n,\alpha}(t,x)k-V^{1,\alpha}(t,x)k)dW_t,\qquad t\in [0,n],\\
e^{-\mu n}(U^{n,\alpha}(n,x)k-U^{2,\alpha}(n,x)k )=e^{-\mu n}U^{n,\alpha}(n,x)k
\end{array}
\right.
\]
We already know that $U^{n,\alpha}(n,x)k$ is uniformly bounded with respect to $n$, 
so the following rate of convergence holds true:
\[
 \vert U^{n,\alpha}(t,x)k-U^{2,\alpha}(t,x)k\vert \leq C e^{-\mu (n-t)}\vert k \vert,
\]
where $C>0$ is a constant that does not depend on $n$. 
So, if we take $(x_j)_{j\geq1},x\in \calh$ such that $x_j\to x$,
then, by the triangular inequality,
\[
\vert U^{2,\alpha}(0,x_j)k-U^{2,\alpha}(0,x)k\vert \leq 2C e^{-\mu (n-t)}\vert k \vert
+\vert U^{n\alpha}(0,x_j)k-U^{n,\alpha}(0,x)k\vert
\]
and, by arguments similar to the ones used in proposition \ref{regback},
the map $x\to U^{n,\alpha}(0,x)k$ is continuous.
So we can conclude that the map $ x\to P^{\alpha}(x)k$ 
is continuous from $\calh$ to
$\mathbb{R}$, linear with respect to $k$ 
and there exists a constant $C>0$ such that
$ |P^{\alpha}(x)k| \leq C |k|_\calh$, so the proof is concluded.
\qed

\begin{corollary}\label{proregu_inf} Setting $v(x)=Y^{x}$, we have
$v\in C(\calh;
\mathbb{R})$ and there exists a constant $C$ such that
$|v(x)|\leq C\, (1+|x|)^2$, $x\in \calh$.
 Moreover $v$ is
G\^{a}teaux differentiable and the map
$(x,h)\to \nabla v(x)h$ is continuous.

For all $\alpha\in [0,1/2)$ and $x\in \calh$ the linear
operator $k\to \nabla v(x)(\lambda-A)^{\alpha}k$ - a priori
defined for $k\in D(\lambda-A)^{\alpha}$ - has an extension to a bounded linear
operator $\calh\to \mathbb{R}$, that we denote by
$[\nabla v(\lambda-A)^{\alpha}](x)$.

Finally the map $(x,k)\to [\nabla v(\lambda-A)^{\alpha}](x)k$
is continuous $ \calh\times \calh \to \mathbb{R}$ and there
exists a constant $C>0$ such that:
\begin{equation}\label{stima_nabla_u_inf}
|[\nabla v (\lambda-A)^{\alpha}](x)k |\leq C |k|_\calh, \qquad x, k\in H.
\end{equation}
\end{corollary}
{\bf Proof.}
We recall that
  $Y^{x}_0$ is
  deterministic. Continuity of $v$ follows from the fact that,
for every $T>0$, the map $x\rightarrow Y_0^{x}$ is
  continuous with values in
  $L^p_{\mathcal{P}}(\Omega,C([0,+\infty],\mathbb{R}))$, $p\geq 2$.

\noindent Similarly, $ \nabla_x v(x)= \nabla_x Y_0^{x}$ exists and
has the required continuity properties, by Proposition
\ref{bsderegolare-inf}.

\noindent Next we notice that
$P^{\alpha}(x)k=
  \nabla_x Y_0^{x}(\lambda-A)^{\alpha}k$. The existence
  of the required extensions
and its continuity are direct consequences of Proposition
\ref{regback_inf}. Finally
the estimate (\ref{stima_nabla_u_inf}) follows from
(\ref{stimadiP_inf}).
\qed
\begin{remark}{\em It is evident by construction
that the law of $Y^{x}$ and consequently the function $v$
depends on the law of the Wiener process $ W$ but not on the
particular probability $\mathbb{P}$ and Wiener process $ W$ we
have chosen. }\end{remark}

\begin{corollary}\label{identmarkov_inf} For every $t\geq 0$, $x\in \calh$
we have, $\P$-a.s.,
\begin{equation}\label{identifY_inf}
Y_s^{x}=v(X_s^{x}),\qquad
{\rm \; for\;all\; }s \geq 0,
\end{equation}
\begin{equation}\label{identifZ_inf}
    Z_s^{x}=
[\nabla v(\lambda-A)^{1-\beta}](X_s^{x})\;(\lambda-A)^{\beta}D_\lambda,\qquad
{\rm \; for\;almost\;all\;}s \geq 0.
\end{equation}
\end{corollary}

{\bf Proof.}
We start from the well-known equality: for $t\geq 0$,
$\P$-a.s.,
$$
X_s^{x}=X_s^{r,X_r^{,x}},\qquad
{\rm \; for\;all\; }s \geq r.
$$
It follows easily from the uniqueness of the backward equation
(\ref{BSDE}) that
$\P$-a.s.,
$$
Y_s^{x}=Y_s ^{r,X_r^{t,x}},\qquad
{\rm \; for\;all\; }s \geq r.
$$
Setting $s=r$ we arrive at (\ref{identifY_inf}).

\noindent To prove (\ref{identifZ_inf}) we note that it follows
immediately from the backward equation
(\ref{BSDE_inf}), see also (\ref{BSDEmild_inf}), that the joint quadratic variation of
$\{Y_s^{x},\;s \geq 0$ and $W$ 
on an arbitrary interval $[t,s]\subset [0,+\infty)$ is equal
to $\int_t^s{Z}_r \; dr$. By
(\ref{identifY_inf}) the same result can be obtained
by considering the joint quadratic variation of
$\{v(X_s^{x}),\;s \geq 0\}$ and $W$. An application
of Proposition \ref{ucompostox} and remark \ref{ucompostox_rem} (whose assumptions hold
true by Corollary \ref{proregu_inf}) leads to the identity
$$
\int_t^s{Z}_r \; dr=
  \int_t^s[\nabla v(\lambda -A)^{1-\beta}](X_r^{x})\,
(\lambda -A)^{\beta}D_\lambda\;dr,
$$
and (\ref{identifZ_inf}) is proved.
\qed

\section{The stationary Hamilton-Jacobi-Bellman equation}

In this section the aim is to solve a second order partial differential equation,
where the second order differential operator is the generator of the Markov process
${X_s^{x},s\geq 0}$, solution of equation (\ref{eqstatonocontrol_inf}).
We denote by $P_{s}$
its transition semigroup:
$$
P_{s}[\phi](x)=\E\, \phi(X_s^x),\qquad x\in \calh,\; s\geq 0,
$$
for any bounded measurable $\phi:\calh\to \R$.
As for the finite horizon case, $P_{s}$ will be considered
as an operator acting on this class of functions.

Let us denote by $\call$ the generator of $P_{s}$,
formally:
$$
\call[\phi](x)=\frac{1}{2}
\< \nabla^2\phi(x)B,B\>
+ \< Ax+F(x),\nabla\phi(x)\>,
$$
where $\nabla\phi(x)$ and $\nabla^2\phi(x)$
are first and second G\^ateaux derivatives of
$\phi$ at the point $x\in \calh$ (here they are
identified with elements of $\calh$ and $L(\calh)$ respectively).

The stationary Hamilton-Jacobi-Bellman equation that we are 
going to study is
\begin{equation}\label{HJBformale_inf}
  \call [v](x) =\mu v(x)-
\Psi (x,\nabla v(t,x)B).
\end{equation}
We consider, for every $T>0$, the variation of
constants formula for (\ref{HJBformale_inf}):
$$
  v(x) =e^{-\mu T}P_{T}[u](x)-\int_0^T e^{-\mu s}P_{s}[
\Psi ( \cdot,
\nabla v(\cdot)B
](x)\; ds,
x\in \calh,
$$
where we recall that $B=(\lambda-A)D_\lambda$.
This equality is still formal, since the
term $(\lambda-A)D_\lambda$ is not defined. However with a slightly
different interpretation we arrive at the following precise
definition:
\begin{definition}\label{defdisoluzionemild_inf}
Let $\beta\in[0,\frac{1}{2})$. We say that a function
$v: \calh\to\R$ is a mild solution of the
Hamilton-Jacobi-Bellman equation
(\ref{HJBformale_inf}) if the following conditions hold:
\begin{enumerate}
  \item[(i)]
$v\in C( \calh;\mathbb{R})$, is
G\^{a}teaux differentiable and the map
$(x,h)\to \nabla v(x)h$ is continuous
$\calh\times \calh\to \R$.

 \item[(iii)]
For all $x\in \calh$ the linear
operator $k\to \nabla v(x)(\lambda-A)^{1-\beta}k$ (a priori
defined for $k\in D(\lambda-A)^{1-\beta}$) has an extension to a bounded linear
operator $\calh\to \mathbb{R}$, that we denote
by $[\nabla v(\lambda-A)^{1-\beta}](x)$.
Moreover the map $(x,k)\to [\nabla v(\lambda-A)^{1-\beta}](x)k$
is continuous $\calh\times \calh \to \mathbb{R}$ and there
exist constants $C,m\ge0$, $\kappa\in [0,1)$ such that
\begin{equation}\label{stima_grad_u_in_generale_inf}
|[\nabla v (\lambda-A)^{1-\beta}](x) |_{\calh^*}\leq C, 
\qquad 
x\in \calh.
\end{equation}

  \item[(iv)] the following equality holds for
  every $x\in \calh$:
  \begin{equation}\label{solmild_inf}
  v(x) =e^{-\mu T}P_{T}[u](x)-\int_0^T e^{-\mu s}P_{s}[
\Psi ( \cdot,
[\nabla v(\lambda-A)^{1-\beta}](\cdot)\;(\lambda-A)^{\beta}D_\lambda )
](x)\; ds,
\end{equation}

\end{enumerate}
\end{definition}


\begin{theorem}\label{main_inf}
Assume Hypotheses \ref{ipotesiconcrete},
\ref{ipsupsi_inf}, \ref{ip_psiagg_inf} and that in equation
(\ref{eqstatonocontrol_inf}) $M$ is taken sufficiently large 
(see also lemma \ref{lemmader_forward_lim}). Then there
exists a unique mild solution of the stationary Hamilton-Jacobi-Bellman
equation (\ref{HJBformale_inf}).
The solution $v$ is given by the formula
$$
    v(x) =Y_0^{x},
  $$
where $(X,Y,Z)$ is the solution of the forward-backward
system (\ref{eqstatonocontrol_inf})-(\ref{BSDEmild_inf}).
\end{theorem}

{\bf Proof.}
The proof is similar to the proof of theorem 6.1 in \cite{FuTe-ell},
noticing, as in \cite{HuTess}, that we can find a mild solution for every
$\lambda>0$, and noticing, as in the finite horizon case, 
see also theorem \ref{main}, that $\nabla v(x) G(x) $ is replaced by
$[\nabla v(\lambda-A)^{1-\beta}](x)\;(\lambda-A)^{\beta}D_\lambda$.

\qed

\section{Synthesis of the optimal control: the infinite horizon case}

Let us consider the cost functional
(\ref{costoconcretoinfor1}), and we make the following assumptions:

\begin{hypothesis}\label{ipotesiconcretecosto_inf}
 $\ell:[0,+\infty)\times\R\times\calu
\to\R$ is continuous and there exists $C>0$ and $g\in L^1 ([0,+\infty))$
such that 
\begin{equation*}
\vert l(\xi,x,u)\vert \leq C g(\xi), \quad \text{for every }\quad
\xi \in [0,+\infty), x\in
\R, u\in\calu.
\end{equation*}
Moreover there exists $C,\epsilon>0$
such that 
\begin{equation*}
\vert l(\xi, x_1,u)-l(\xi, x_2,u)\vert\leq 
C \frac{\vert x_1-x_2\vert}{(1+\xi)^{\frac{1+\epsilon}{2}}}\sqrt{\rho(\xi)}
\quad \text{for every }\quad
\xi \in [0,+\infty), x_1,x_2\in\R, u\in\calu.
\end{equation*}

\end{hypothesis}

In this section we assume that Hypothesis \ref{ipotesiconcretecosto_inf}
holds.
We briefly reformulate the cost (\ref{costoconcretoinfor1}) in an abstract form.
We define
$$
L(x,u)=\int_0^{+\infty} \ell(s,\xi, x(\xi),u)\;d\xi, \qquad
 x=x(\cdot)\in \calh,u\in \calu,$$
and so $L: \calh
\times \calu\to\R$ is well defined and measurable and the cost functional 
(\ref{costoconcretoinfor1}) can be written in the form
\begin{equation}\label{costoastratto_inf}
J(x,u(\cdot))=\E \int_0^{+\infty}e^{-\mu s} L( X_s^u,u_s)\;ds.
\end{equation}
Moreover for $x\in \calh$, $z\in\R$ we define the
hamiltonian:
$$
\Psi(x,z)=\inf_{u\in\calu} \{zu+L(x ,u) \},
$$
where $zu$ denotes the scalar product in $\R$. 
We notice that setting $C_\calu=\sup\{|u|\,:\,u\in\calu\}$ we have
  $|\Psi(x,z_1)-\Psi(x,z_2)|\le C_\calu\, |z_1-z_2|, $ for every
 $x\in \calh$, $z_1,z_2\in\R$.
Moreover we assume that the hamiltonian $\Psi$ satisfies
hypothesis \ref{ip_psiagg_inf}.

\noindent Analougsly to the infinite horizon case, if we define
\begin{equation}\label{defdigammagrande_inf}
\Gamma(x,z)=\left\{ u\in \calu: zu+L(x ,u)= \Psi(x,z)\right\}
\end{equation}
 then
$\Gamma(x,z) \neq \emptyset$ for every $x\in
\calh$ and every $z\in \R$ and so it admits a measurable selection,
$\gamma: \calh \times \R 
\rightarrow \calu$ with $\gamma$  measurable and $\gamma(x,z)\in \Gamma(x,z)$ 
for every $x\in \calh$ and every $z\in \R$.

We now reformulate the optimal control problem in the
weak sense, following the approach
of \cite{FlSo}. As in section 6.3,
$(\Omega,\calf, ({\cal F}_{t}), \P, W)$ is
an {\it admissible set-up}, and  
$\U=(\Omega,\calf, ({\cal F}_{t}), \P, W, u, X^u)$ is an
{\it admissible control
system} (a.c.s.) if:
 \begin{itemize}
    \item $(\Omega,\calf, ({\cal F}_{t}), \P, W)$
is an admissible set-up;

    \item $u:\Omega\times [0,+\infty)\to\R$ is an
$({\cal F}_{t})$-predictable process
  with values in $\calu$;

\item $\{X^u_t,\, t\in [0,+\infty)\}$ is an
$({\cal F}_{t})$-adapted continuous process
  with values in $\calh$, mild solution
of the state equation (\ref{eqstatoformale_inf})
 with initial condition
 $X^{u}_0=x$.
 \end{itemize}

For $x\in \calh$ we wish to minimize
the cost (\ref{costoastratto_inf}):
\begin{equation}\label{costoastrattodue_inf}
J(x,U(\cdot))=\E \int_0^{+\infty} e^{-\mu s}L(X_s^u,u(s))\;ds 
\end{equation}
over all
admissible control systems.
%

We recall that by $v:\calh\to\R$, we denote the mild solution of
the Hamilton-Jacobi-Bellman equation (\ref{HJBformale_inf}).

\begin{theorem}\label{th-rel-font-inf}
 Assume Hypotheses \ref{ipotesiconcrete},
\ref{ipotesiconcretecosto_inf} and hat $\Psi$
satisfies hypothesis \ref{ip_psiagg_inf}.
For every $x\in \calh$ and for
 all admissible control systems $U$ we have $J(t,x,U))
 \geq v(x)$,
  and the
 equality holds if and only if
$$
u_s\in \Gamma\left(,X^{u}_s, [\nabla
v(\lambda-A)^{1-\beta}](s ,X^{u}_s)
\;(\lambda-A)^{\beta}D_\lambda\right)
  $$
Moreover, if
\begin{equation*}
u(x)=\gamma\Big(x, [\nabla
v(\lambda-A)^{1-\beta}](x)\;(\lambda-A)^{\beta}D_\lambda \Big),\qquad
x\in \calh,
\end{equation*}
then there exists an adapted process $\{\overline{X}_s,\;s\geq 0\}$ with 
continuous trajectories solving the
closed loop equation: $\P$-a.s.
\begin{equation}\label{cle_inf}
\begin{array}{lll}\dis
\overline{X}_s&=&\dis e^{sA}x_0+\int_{0}^s
e^{(s-r)A}F(r,\overline{X}_r)\;dr
+
\int_t^s(\lambda -A)^{1-\beta}e^{(s-r)A}(\lambda
-A)^{\beta}D_\lambda\;dW_r
\\&
+&\dis\int_{0}^s(\lambda -A)^{1-\beta}e^{(s-r)A}(\lambda
-A)^{\beta}D_\lambda\;u(r,\overline{X}_r)dr,
\qquad s\geq 0,
\end{array}
\end{equation}
and $(\overline{X}_s, \gamma(\overline{X}_s,[\nabla
v(\lambda-A)^{1-\beta}](\overline{X}_s)\;(\lambda-A)^{\beta}D_\lambda ))$
is an optimal pair.
\end{theorem}

{\bf Proof.} The proof is similar to the proof of Theorem 5.1 in \cite{HuTess}. Just notice that in this case by
(\ref{identifZ_inf}) we have $Z_s^{x}= [\nabla
v(\lambda-A)^{1-\beta}](X_s ^{x})\;(\lambda-A)^{\beta}D_\lambda$ and
the role of $G$ in \cite{HuTess}, Theorem 5.1 is here played by
$B=(\lambda-A)D_\lambda$. \qed

\begin{remark}
 We notice that the techniques used to treat the stationary Hamilton Jacobi Bellman 
equation and the infinite horizon optimal control problem can be applied to the case of boundary conditions of 
Neumann type in the state equation, i.e to a state equation like the one studied in \cite{DebFuTe} but considered 
for every $t>0$.
\end{remark}


\end{document}